%% file: Sp.tex
\documentclass[11pt,reqno]{amsart}
\usepackage{amssymb,amscd,amsbsy}
\usepackage{amssymb,amscd,amsbsy,mathrsfs}
\input{classstartscr}

\theoremstyle{remark}

\newtheorem*{rem*}{Remark}

\newcommand\up{\upsilon}
\newcommand\Li{{\rm Lip}}
\newcommand\fM{\frak M}
\newcommand\cZ{\mathcal{Z}}
\newcommand\dg{\frak D}

\begin{document}

\newcommand{\qm}{\quad\mbox{and}\quad}
\newcommand{\vse}{\vspace{.2in}}
\numberwithin{equation}{section}

\title{Functions of operators under perturbations of class $\bS_p$}
\author{A.B. Aleksandrov and V.V. Peller}
\thanks{The first author is partially supported by RFBR grant 08-01-00358-a and by
Russian Federation presidential grant NSh-2409.2008.1;
the second author is partially supported by NSF grant DMS 0700995 and by ARC grant}
\maketitle

\

\begin{center}
{\Large Contents}
\end{center}

\

\begin{abstract}
This is a continuation of our paper \cite{AP2}. We prove that for functions $f$ in the H\"older class $\L_\a(\R)$ and $1<p<\be$, the operator 
$f(A)-f(B)$ belongs to $\bS_{p/\a}$, whenever $A$ and $B$ are self-adjoint operators with $A-B\in\bS_p$. We also obtain sharp estimates for the Schatten--von Neumann norms $\big\|f(A)-f(B)\big\|_{\bS_{p/\a}}$ in terms of $\|A-B\|_{\bS_p}$ and establish similar results for other operator ideals. We also estimate Schatten--von Neumann norms of higher order differences 
$\sum\limits_{j=0}^m(-1)^{m-j}\left(\begin{matrix}m\\j\end{matrix}\right)f\big(A+jK\big)$. We prove that analogous results hold for functions on the unit circle and unitary operators and for analytic functions in the unit disk and contractions.
Then we find necessary conditions on $f$ for $f(A)-f(B)$ to belong to $\bS_q$ under the assumption that $A-B\in\bS_p$. We also obtain Schatten--von Neumann estimates for quasicommutators $f(A)Q-Qf(B)$, and introduce a spectral shift function and find a trace formula for operators of the form 
$f(A-K)-2f(A)+f(A+K)$.
\end{abstract}

\

\begin{enumerate}
\item[1.] Introduction \quad\dotfill \pageref{In}
\item[2.] Besov spaces \quad\dotfill \pageref{fs}
\item[3.] Ideals of operators on Hilbert space \quad\dotfill \pageref{ide}  
\item[4.] Multiple operator integrals \quad\dotfill \pageref{koi}
\item[5.] Self-adjoint operators. Sufficient conditions \quad\dotfill \pageref{sa}
\item[6.]  Unitary operators. Sufficient conditions \quad\dotfill \pageref{u}
\item[7.]  The case of contractions \quad\dotfill \pageref{con}
\item[8.] Finite rank perturbations and necessary conditions. Unitary operators
\quad\dotfill \pageref{fru}
\item[9.] Finite rank perturbations and necessary conditions. Self-adjoint operators \quad\dotfill \pageref{frsa}
\item[10.]  Spectral shift function for second order differences\quad\dotfill \pageref{fss}
\item[11.]  Commutators and quasicommutators\quad\dotfill \pageref{cqc}
\item[] References \quad\dotfill \pageref{bibl}
\end{enumerate}

\

\newcommand{\mt}{{\mathscr T}}
\newcommand{\fI}{{\frak I}}

\setcounter{section}{0}
\section{\bf Introduction}
\setcounter{equation}{0}
\label{In}

\

This paper is a continuation of our paper \cite{AP2}. In \cite{AP2} we obtained sharp estimates for the norms of $f(A)-f(B)$ in terms of the norm of $A-B$ for various classes of functions $f$. Here $A$ and $B$ are self-adjoint operators on Hilbert space and $f$ is a function on the real line $\R$. We also obtained in \cite{AP2} sharp estimates for the norms of higher order differences
\bay
\label{hod}
\big(\D_K^mf\big)(A)\df
\sum_{j=0}^m(-1)^{m-j}\left(\begin{matrix}m\\j\end{matrix}\right)f\big(A+jK\big),
\ey
where $A$ and $K$ are self-adjoint operators. Similar results were obtained in \cite{AP2} for functions of unitary operators and for functions of contractions.

In this section we are going to obtain sharp estimates for the Schatten--von Neumann norms of first order differences $f(A)-f(B)$ and higher order differences
$\big(\D_K^mf\big)(A)$ for functions $f$ that belong to a H\"older--Zygmund class
$\L_\a(\R)$, $0<\a<\be$, (see \S\,\ref{fs} for the definition of these spaces).

In particular we study the question, under which conditions on $f$ the operator
\lb$f(A)-f(B)$ (or $\big(\D_K^mf\big)(A)\,$) belongs to the Schatten--von Neumann class $\bS_q$, whenever $A-B$ (or $K$) belongs to $\bS_p$.

We also obtain related results for more general ideals of operators on Hilbert space (see \S\,\ref{ide} for the introduction to operator ideals on Hilbert space). 

In connection with the Lifshits--Krein trace formula, M.G.~Krein asked in 
\cite{Kr} \lb the question whether $f(A)-f(B)\in\bS_1$, whenever $f$ is a {\it Lipschitz function} (i.e., \lb$|f(x)-f(y)|\le\const|x-y|$, $x,\,y\in\R$) 
and $A-B\in\bS_1$. Functions $f$ satisfying this property are called {\it trace class perturbations preserving}.

Farforovskaya constructed in \cite{F} an example that shows that the answer to the Krein question is negative. 

Later in \cite{Pe2} and \cite{Pe4} necessary conditions and sufficient conditions
for $f$ to be trace class perturbations preserving were found. It was shown in
\cite{Pe2} and \cite{Pe4} that if $f$ belongs to the Besov space $B^1_{\be1}(\R)$ (see \S\,\ref{fs}), then $f$ is trace class perturbations preserving. On the other hand, it was shown in \cite{Pe2} that if $f$ is trace class perturbations preserving, then it belongs to the Besov space $B_1^1(\R)$ locally.
This necessary condition also proves that a Lipschitz function does not have to be trace class perturbations preserving. Moreover, in \cite{Pe2} and \cite{Pe4}
a stronger necessary condition was also found. Note that a function is trace class perturbations preserving if and only if it is operator Lipschitz
(see \cite{Pe2} and \cite{KS}).

We also mention here the paper \cite{Pe3}, in which analogs of the above results 
were obtained for perturbations of class $\bS_p$ with $p\in(0,1)$.

On the other hand, Birman and Solomyak in \cite{BS3} proved that a Lipschitz function $f$ must preserve Hilbert--Schimidt class perturbations: $f(A)-f(B)\in\bS_2$, whenever $A-B\in\bS_2$ and
$$
\|f(A)-f(B)\|_{\bS_2}\le\sup_{x\ne y}\frac{|f(x)-f(y)|}{|x-y|}\|A-B\|_{\bS_2}.
$$
To prove that result, Birman and Solomyak developed in \cite{BS1}, \cite{BS2}, and \cite{BS3} their beautiful theory of double operator integrals and established a formula for $f(A)-f(B)$ in terms of double operator integrals
(see \S\,\ref{koi}). Note also that the paper \cite{KS} studies functions that preserve perturbations belonging to operator ideals.

We mention here two recent results. In \cite{NP} it was proved that if $f$ is a Lipschitz function and $\rank(A-B)<\be$, then $f(A)-f(B)$ belongs the weak space $\bS_{1,\be}$ (see \S\,\ref{ide} for the definition). It was also shown in \cite{NP} that if $A-B\in\bS_1$, then $f(A)-f(B)$ belongs to
the ideal $\bS_\O$, i.e.,
$$
\sum_{j=0}^ns_j\big(f(A)-f(B)\big)\le\const\log(2+n).
$$
(here $s_j$ is the $j$th singular value).
This allowed the authors of \cite{NP} to deduce that for $p\ge1$ and $\e>0$,
the operator $f(A)-f(B)$ belongs to $\bS_{p+\e}$, whenever $f$ is a Lipschitz function and $A-B\in\bS_p$.

The epsilon was removed later in \cite{PS} in the case $1<p<\be$. It was shown in \cite{PS} that for $p\in(1,\be)$, the operator $f(A)-f(B)$ belongs to $\bS_p$, whenever $A-B\in\bS_p$ and $f$ is a Lipschitz function.

Note that similar results also hold for functions on the unit circle $\T$ and unitary operators.

It was shown in \cite{BKS}  that if $A$ and $B$ are positive self-adjoint operators and $\fI$ is a normed ideal of operators on Hilbert space with majorization property, then for $\a\in(0,1)$, the following inequality holds:
$$
\big\|A^\a-B^\a\big\|_\fI\le\big\|\,|A-B|^\a\big\|_\fI.
$$
In this paper we study the problem under which conditions on a function
$f$ and a (quasi)normed ideal $\fI$ of operators on Hilbert space the following
inequality holds:
$$
\big\|f(A)-f(B)\big\|_\fI\le\const\big\|\,|A-B|^\a\big\|_\fI.
$$

In Section \ref{sa} of this paper among other results we show that if $f$
belongs to the H\"older class $\L_\a$, $0<\a<1$, and $1<p<\be$,
then $f(A)-f(B)\in\bS_{p/\a}$ and
$$
\big\|f(A)-f(B)\big\|_{\bS_{p/\a}}\le\const\|f\|_{\L_\a(\R)}\|A-B\|_{\bS_p}^\a.
$$
On the other hand, this is not true for $p=1$ (a counter-example is given in 
\S\,\ref{frsa}). Nevertheless, for $p=1$, under the assumptions that 
$f\in\L_\a(\R)$ and $A-B\in\bS_1$, we prove that $f(A)-f(B)$ belongs to the weak space $\bS_{\frac{1}{\a},\be}$. To make the conclusion that 
$f(B)-f(A)\in\bS_{1/\a}$ under the assumption that $A-B\in\bS_1$, we need the stronger condition: $f$ belongs to the Besov space $B_{\be1}^\a$. We also obtain similar results for other ideals of operators on Hilbert space. In particular,
we show that for every $p\in(1,\be)$ and every $l\ge0$, the following inequality holds
$$
\sum_{j=0}^l\Big(s_j\big(|f(A)-f(B)|^{1/\a}\big)\Big)^p
\le\const\|f\|^{p/\a}_{\L_\a(\R)}\sum_{j=0}^l\big(s_j(A-B)\big)^p,
$$
where the constant does not depend on $l$. We also establish in \S\,\ref{sa} similar results for higher order differences
$\big(\D_K^mf\big)(A)$ and functions $f\in\L_\a(\R)$ with $\a\in[m-1,m)$.

In \S\,\ref{u} we obtain analogs of the result of \S\,\ref{sa} for functions on $\T$ and unitary operators, while in \S\,\ref{con} we establish similar results for functions analytic in the unit disk and contractions.

In Section \ref{fru} we obtain refinements of some results of \S\,\ref{u} in the case of finite rank perturbations of unitary operators. We also give some necessary conditions on a function $f$ for $f(U)-f(V)$ to belong to $\bS_q$, whenever $U-V\in\bS_p$. Analogs of the results of \S\,\ref{fru} for self-adjoint operators are given in \S\,\ref{frsa}.

In \S\,\ref{fss} we consider the problem of evaluating the trace of $f(A-K)-2f(A)+f(A+K)$ under the assumptions that $K\in\bS_2$ and $f$ belongs
to the Besov class $B_{\be1}^2(\R)$. We introduce a spectral shift function
$\varsigma$ associated with the pair $(A,K)$ and establish the following trace formula:
$$
\trace\big(f(A-K)-2f(A)+f(A+K)\big)=\int_\R f''(x)\varsigma(x)\,dx.
$$
We also show that similar results hold in the case of unitary operators.

The final section \ref{cqc} is devoted to estimates of commutators and quasicommutators
in the norm of Schatten--von Neumann classes (as well as in the norms of more general operator ideals). We consider a bounded operator $Q$, self-adjoint operators $A$ and $B$ and for a function $f\in\L_\a(\R)$, we prove that $f(A)Q-Qf(B)\in\bS_{p/\a}$, whenever $p>1$ and $AQ-QB\in\bS_p$. We also obtain
norm estimates for $f(A)Q-Qf(B)$ that are similar to the estimates obtained in 
\S\,\ref{sa} for first order differences $f(A)-f(B)$. 


In \S\,\ref{fs} we give a brief introduction to Besov spaces and, in particular, we discuss H\"older--Zygmund classes $\L_\a(\R)$, $0<\a<\be$.

In Section \ref{ide} we introduce quasinormed ideals of operators on Hilbert space and define the upper Boyd index of a quasinormed ideal. 

In \S\,\ref{koi} we give an introduction to the Birman--Solomyak theory of double operator integrals which will be used in the paper to obtain desired estimates. We also define multiple operator integrals and multiple operator integrals with respect to semi-spectral measures. 

Note that in this paper we give detailed proofs in the case of bounded self-adjoint operators and explain briefly that the same results also hold in the case of unbounded self-adjoint operators. We are going to consider in detail the case of unbounded self-adjoint operators in \cite{AP3}. Note also that we are going to
consider separately in \cite{AP4} similar problems for perturbations of dissipative operators and improve earlier results of \cite{Nab}.

The main results of this paper have been announced without proofs in \cite{AP1}.

\

\section{\bf Besov spaces}
\setcounter{equation}{0}
\label{fs}

\

The purpose of this section is to give a brief introduction to Besov spaces that play an important role in problems of perturbation theory.
We start with Besov spaces on the unit circle.

Let $1\le p,\,q\le\be$ and $s\in\R$. The Besov class $B^s_{pq}$ of functions (or
distributions) on $\T$ can be defined in the following way. Let $w$ be an infinitely differentiable function on $\R$ such
that
\bay
\label{w}
w\ge0,\quad\supp w\subset\left[\frac12,2\right],\quad\mbox{and} \quad w(x)=1-w\left(\frac x2\right)\quad\mbox{for}\quad x\in[1,2].
\ey

Consider the trigonometric polynomials $W_n$, and $W_n^\sharp$ defined by
$$
W_n(z)=\sum_{k\in\Z}w\left(\frac{k}{2^n}\right)z^k,\quad n\ge1,\quad W_0(z)=\bar z+1+z,\quad
\mbox{and}\quad W_n^\sharp(z)=\ov{W_n(z)},\quad n\ge0.
$$
Then for each distribution $f$ on $\T$,
$$
f=\sum_{n\ge0}f*W_n+\sum_{n\ge1}f*W^\sharp_n.
$$
The Besov class $B^s_{pq}$ consists of functions (in the case $s>0$) or distributions $f$ on $\T$
such that
\bay
\label{bes}
\big\{\|2^{ns}f*W_n\|_{L^p}\big\}_{n\ge1}\in\ell^q\quad\mbox{and}
\quad\big\{\|2^{ns}f*W^\sharp_n\|_{L^p}\big\}_{n\ge1}\in\ell^q.
\ey


Besov classes admit many other descriptions. In particular, for $s>0$, the space $B^s_{pq}$ admits the
following characterization. A function $f\in L^p$ belongs to $B^s_{pq}$, $s>0$, if and only if
$$
\int_\T\frac{\|\D^n_\t f\|_{L^p}^q}{|1-\t|^{1+sq}}d\m(\t)<\be\quad\mbox{for}\quad q<\be
$$
and
\bay
\label{pbe}
\sup_{\t\ne1}\frac{\|\D^n_\t f\|_{L^p}}{|1-\t|^s}<\be\quad\mbox{for}\quad q=\be,
\ey
where $\m$ is normalized Lebesgue measure on $\T$, $n$ is an integer greater than $s$, and $\D_\t$, $\t\in\T$, is
the difference operator: 
$$
(\D_\t f)(\z)=f(\t\z)-f(\z), \quad\z\in\T.
$$

We use the notation $B_p^s$ for $B_{pp}^s$.

The spaces $\L_\a\df B_\be^\a$ form the {\it H\"older--Zygmund scale}. If $0<\a<1$, then $f\in\L_\a$ if and only if
$$
|f(\z)-f(\t)|\le\const|\z-\t|^\a,\quad\z,\,\t\in\T,
$$
while $f\in\L_1$ if and only if $f$ is continuous and
$$
|f(\z\t)-2f(\z)+f(\z\bar\t)|\le\const|1-\t|,\quad\z,\,\t\in\T.
$$
By \rf{pbe},  $\a>0$, $f\in\L_\a$ if and only if $f$ is continuous and
$$
|(\D^n_\t f)(\z)|\le\const|1-\t|^\a,
$$
where $n$ is a positive integer such that $n>\a$. 

Note that the (semi)norm of a function $f$ in $\L_\a$ is equivalent to
$$
\sup_{n\ge1}2^{n\a}\big(\|f*W_n\|_{L^\be}+\|f*W_n^\sharp\|_{L^\be}\big).
$$



It is easy to see from the definition of Besov classes that the Riesz projection $\pp_+$,
$$
\pp_+f=\sum_{n\ge0}\hat f(n)z^n,
$$
is bounded on $B^s_{pq}$. Functions in $\big(B^s_{pq}\big)_+\df\pp_+B^s_{pq}$ admit a natural extension to analytic functions
in the unit disk $\dd$. It is well known that the functions in $\big(B^s_{pq}\big)_+$ admit the following description:
$$
f\in \big(B^s_{pq}\big)_+\Leftrightarrow
\int_0^1(1-r)^{q(n-s)-1}\|f^{(n)}_r\|^q_p\,dr<\be,\quad q<\be,
$$
and
$$
f\in \big(B^s_{p\be}\big)_+\Leftrightarrow
\sup_{0<r<1}(1-r)^{n-s}\|f^{(n)}_r\|_p<\be,
$$
where $f_r(\z)\df f(r\z)$ and $n$ is a nonnegative integer greater than $s$.

Let us proceed now to Besov spaces on the real line. We consider homogeneous Besov spaces 
$B_{pq}^s(\R)$ of functions (distributions) on $\R$.
We use the same function $w$
as in \rf{w} and define the functions $W_n$ and $W^\sharp_n$ on $\R$ by
$$
\F W_n(x)=w\left(\frac{x}{2^n}\right),\quad\F W^\sharp_n(x)=\F W_n(-x),\quad n\in\Z,
$$
where $\F$ is the {\it Fourier transform}:
$$
\big(\F f\big)(t)=\int_\R f(x)e^{-{\rm i}xt}\,dx,\quad f\in L^1.
$$

With every tempered distribution $f\in{\mathscr S}^\prime(\R)$ we 
associate a sequences $\{f_n\}_{n\in\Z}$,
$$
f_n\df f*W_n+f*W_n^\sharp.
$$
Initially we define the (homogeneous) Besov class $\dot B^s_{pq}(\R)$ as the set of all $f\in{\mathscr S}^\prime(\R)$
such that 
\bay
\label{Wn}
\{2^{ns}\|f_n\|_{L^p}\}_{n\in\Z}\in\ell^q(\Z).
\ey
According to this definition, the space $\dot B^s_{pq}(\R)$ contains all polynomials. Moreover, the distribution $f$ is defined by the sequence $\{f_n\}_{n\in\Z}$
uniquely up to a polynomial. It is easy to see that the series $\sum_{n\ge0}f_n$ converges in ${\mathscr S}^\prime(\R)$.
However, the series $\sum_{n<0}f_n$ can diverge in general. It is easy to prove that the
series $\sum_{n<0}f_n^{(r)}$ converges uniformly on $\R$ for each nonnegative integer
$r>s-1/p$. Note that in the case $q=1$ the series $\sum_{n<0}f_n^{(r)}$ converges uniformly whenever $r\ge s-1/p$.

Now we can define the modified (homogeneous) Besov class $B^s_{pq}(\R)$. We say that a distribution $f$
belongs to $B^s_{pq}(\R)$ if $\{2^{ns}\|f_n\|_{L^p}\}_{n\in\Z}\in\ell^q(\Z)$ and
$f^{(r)}=\sum_{n\in\Z}f_n^{(r)}$ in the space ${\mathscr S}^\prime(\R)$, where $r$ is
the minimal nonnegative integer such that $r>s-1/p$ ($r\ge s-1/p$ if $q=1$). Now the function $f$ is determined uniquely by the sequence $\{f_n\}_{n\in\Z}$ up
to a polynomial of degree less that $r$, and a polynomial $\varphi$ belongs to $B^s_{pq}(\R)$
if and only if $\deg\varphi<r$.


We use the same notation $W_n$ and $W_n^\sharp$ for functions on $\T$ and on $\R$. This will not lead to a confusion.

Besov spaces $B^s_{pq}(\R)$ admit equivalent definitions that are similar to those discussed above in the case of Besov
spaces of functions on $\T$. In particular, the H\"older--Zygmund classes $\L_\a(\R)\df B^\a_{\be}(\R)$, $\a>0$, can be described 
as the classes of continuous functions $f$ on $\R$ such that
$$
\big|(\D^m_tf)(x)\big|\le\const|t|^\a,\quad t\in\R,
$$
where the difference operator $\D_t$ is defined by
$$
(\D_tf)(x)=f(x+t)-f(x),\quad x\in\R,
$$
and $m$ is an integer greater than $\a$.

As in the case of functions on the unit circle, we consider the following (semi)norm on $\L_\a(\R)$:
$$
\sup_{n\in\Z}2^{n\a}\big(\|f*W_n\|_{L^\be}+\|f*W_n^\sharp\|_{L^\be}\big),\quad f\in\L_\a(\R).
$$

We refer the reader to \cite{Pee} and \cite{Pe5} for more detailed information on Besov spaces.

\

\section{\bf Ideals of operators on Hilbert space}
\setcounter{equation}{0}
\label{ide}

\

In this section we give a brief introduction to quasinormed ideals of operators on Hilbert space.
First we recall the definition of quasinormed vector spaces.

Let $X$ be a vector space. A functional $\|\cdot\|:X\to[0,\be)$ is called a {\it quasinorm} on $X$
if 

(i) $\|x\|=0$ if and only if $x=\0$;

(ii) $\|\a x\|=|\a|\cdot\|x\|$, for every $x\in X$ and $\a\in\C$;

(iii) there exists a positive number $c$ such that $\|x+y\|\le c\big(\|x\|+\|y\|)$ for every $x$ and $y$ in 
$X$.

We say that a sequence $\{x_j\}_{j\ge1}$ of vectors of a {\it quasinormed space} $X$ converges to 
$x\in X$ if $\lim\limits_{j\to\be}\|x_j-x\|=0$. It is well known that there exists a translation invariant metric on $X$ which induces an equivalent topology on $X$. A quasinormed space is called {\it quasi-Banach} if it is complete.

To proceed to operator ideals on Hilbert space, we also recall the definition of singular values of bounded linear operators on Hilbert space. Let $T$ be a bounded linear operator. The singular
values $s_j(T)$, $j\ge0$, are defined by
$$
s_j(T)=\inf\big\{\|T-R\|:~\rank R\le j\big\}.
$$
Clearly, $s_0(T)=\|T\|$ and $T$ is compact if and only if $s_j(T)\to0$ as $j\to\be$.

For a bounded operator $T$ on Hilbert space we also introduce the sequence $\{\s_n(T)\}_{n\ge0}$
defined by
\bay
\label{sin}
\s_n(T)\df\frac1{n+1}\sum_{j=0}^ns_j(T).
\ey

\medskip

{\bf Definition.} Let $\h$ be a Hilbert space and let ${\frak I}$ be a linear manifold in the set $\B(\h)$ of bounded linear operators on $\h$ that is equipped with a quasi-norm $\|\cdot\|_{\frak I}$ that makes
$\fI$ a quasi-Banach space.
We say that ${\frak I}$ is a {\it quasinormed ideal} if for every $A$ and $B$ in $\B(\h)$ and 
$T\in{\frak I}$,
\bay
\label{qni}
ATB\in{\frak I}\qm\|ATB\|_\fI\le\|A\|\cdot\|B\|\cdot\|T\|_\fI.
\ey
A quasinormed ideal $\fI$ is called a {\it normed ideal} if $\|\cdot\|_\fI$ is a norm.

Note that we do not require that $\fI\ne\B(\h)$.

\medskip

It is easy to see that if $T_1$ and $T_2$ are operators in a quasinormed ideal $\fI$ and 
$s_j(T_1)=s_j(T_2)$ for $j\ge0$, then $\|T_1\|_\fI=\|T_2\|_\fI$. Thus there exists a function 
$\Psi=\Psi_\fI$ defined on the set of nonincreasing sequences of nonnegative real numbers with values in $[0,\be]$
such that $T\in\fI$ if and only if $\Psi\big(s_0(T),s_1(T),s_2(T),\cdots~)<\be$ and
$$
\|T\|_\fI=\Psi\big(s_0(T),s_1(T),s_2(T),\cdots~),\quad T\in\fI.
$$
If $T$ is an operator from a Hilbert space $\h_1$ to a Hilbert space $\h_2$, we say that $T$ belongs to 
$\fI$ if $\Psi\big(s_0(T),s_1(T),s_2(T),\cdots~)<\be$.

For a quasinormed ideal $\fI$ and a positive number $p$, we define the quasinormed ideal $\fI^{\{p\}}$ by
$$
\fI^{\{p\}}=\left\{T:~\big(T^*T\big)^{p/2}\in\fI\right\},
\quad\|T\|_{\fI^{\{p\}}}\df\left\|(T^*T\big)^{p/2}\right\|_\fI^{1/p}.
$$

If $T$ is an operator on a Hilbert space $\h$ and $d$ is a positive integer, we denote by $[T]_d$
the operator on $\bigoplus\limits_{j=1}^dT_j$ on the orthogonal sum of $d$ copies of $\h$, where 
$T_j=T$, $1\le j\le d$. It is easy to see that 
$$
s_n\big([T]_d\big)=s_{[n/d]}(T),\quad n\ge0,
$$
where $[x]$ denotes the largest integer that is less than or equal to $x$.

We denote by $\b_{\fI,d}$ the quasinorm of the transformer $T\mapsto[T]_d$
on $\fI$. Clearly, the sequence $\{\b_{\fI,d}\}_{d\ge1}$ is nondecreasing and {\it submultiplicative}, i.e.,
$\b_{\fI,d_1d_2}\le\b_{\fI,d_1}\b_{\fI,d_2}$. It is well known that the last inequality implies that
\bay
\label{lif}
\lim_{d\to\be}\frac{\log\b_{\fI,d}}{\log d}=\inf_{d\ge2}\frac{\log\b_{\fI,d}}{\log d}.
\ey

An analog of \rf{lif} for submultiplicative functions on $(0,\be)$ is proved in \cite{KPS}, Ch. II, Th. 1.3.
To reduce the case of sequences to the case of functions, one can proceed as follows.
Suppose that $\{\b_n\}_{n\ge1}$ is a nondecreasing submultiplicative sequence such that $\b_1=1$.
We can define the function $v$ on $(0,\be)$ by $v(t)=\min\{\b_n:~n\ge t\}$. Then $v(n)=\b_n$ and to prove \rf{lif}, it suffices to apply Theorem 1.3 of Ch. 2 of \cite{KPS} to the function $v$.

\medskip

{\bf Definition.} If $\fI$ is a quasinormed ideal, the number
$$
\b_\fI\df\lim_{d\to\be}\frac{\log\b_{\fI,d}}{\log d}=\inf_{d\ge2}\frac{\log\b_{\fI,d}}{\log d}
$$
is called the {\it upper Boyd index of} $\fI$.

\medskip

It is easy to see that $\b_\fI\le1$ for an arbitrary normed ideal $\fI$. It is also clear that $\b_\fI<1$
if and only if $\lim\limits_{d\to\be}d^{-1}\b_{\fI,d}=0$.

Note that the upper Boyd index does not change if we replace the initial quasinorm in the quasinormed ideal with an equivalent one that also satisfies \rf{qni}. It is also easy to see that
$$
\b_{\fI^{\{p\}}}=p^{-1}\b_\fI.
$$

Theorem \ref{Bo} below is known to experts. Its analog for symmetrically normed spaces can be found in \cite{KPS}, Ch. 2, Th. 6.6. A similar method can be used to prove Theorem \ref{Bo}. We give a proof here for reader's convenience.

\begin{thm}
\label{Bo}
Let $\fI$ be a quasinormed ideal. The following are equivalent:

{\em(i)} $\b_\fI<1$;

{\em(ii)} for every nonincreasing sequence $\{s_n\}_{\ge0}$
of nonnegative numbers,
\bay
\label{xi}
\Psi_\fI\Big(\{\s_n\}_{n\ge0}\Big)\le\const\Psi_\fI\Big(\{s_n\}_{n\ge0}\Big),
\ey
where $\s_n\df(1+n)^{-1}\sum\limits_{j=0}^ns_j$.
\end{thm}

In the proof of Theorem \rf{Bo} we are going to use an elementary fact that if $\sum\limits_{n\ge1}x_n$
is a series of vectors in a quasi-Banach space $X$ such that $\|x_n\|\le\const\g^n$ for some $\g<1$, then the series converges in $X$. This is obvious if $c\g<1$, where $c$ is the constant in the definition of quasinorms. In the general case we can partition the series $\sum\limits_{n\ge1}x_n$ in several 
series, after which each resulting series satisfies the above assumption.

\medskip

{\bf Proof of Theorem \ref{Bo}.}  Let us first show that (i)$\imp$(ii). Suppose that $\b_\fI<\d<1$. Then there exists $C>0$ such that
$\b_{\fI,d}\le C\d^d$ for all positive $d$. Let $\{s_n\}_{n\ge0}$ be a nonincreasing sequence of positive numbers such that $\Psi\big(\{s_n\}_{n\ge0}\big)<\be$. Let $\{e_j\}_{j\ge0}$ be an orthonormal basis in a Hilbert space $\h$. For $k\ge0$, we consider the operator $A_k\in\B(\h)$ defined by
$A_ke_j=s_{[2^{-k}j]}e_j$, $j\ge0$. It is easy to see that $A_k$ is unitarily equivalent to the operator
$\big[A_0\big]_{2^k}$ and
$$
\|A_k\|_\fI\le C2^{\d k}\|A_0\|_\fI=C2^{\d k}\Psi\big(\{s_n\}_{n\ge0}\big).
$$
It follows that the series $A=\sum\limits_{k=0}^\be2^{-k}A_k$ converges in $\fI$ and $\|A\|_\fI\le c\,\Psi\big(\{s_n\}_{n\ge0}\big)$,
where $c$ is a positive number. Clearly,
$s_n(A)=\sum\limits_{k=0}^\infty 2^{-k}s_{\big[\frac n{2^k}\big]}$.

We have
\begin{align*}
\sum_{j=0}^ns_j&=s_n+\sum_{k=1}^{1+[\log_2n]}~\sum_{j=[2^{-k}n]}^{[2^{-k+1}n]-1}s_j
\le s_n+\sum_{k=1}^{1+[\log_2n]}\big([2^{-k+1}n]-[2^{-k}n]\big)s_{[2^{-k}n]}\\[.2cm]
&\le s_n+\sum_{k=1}^{1+[\log_2n]}\big(2^{-k}n+1\big)s_{[2^{-k}n]}
\le s_n+3n\sum_{k=1}^\be 2^{-k}s_{[2^{-k}n]}\\[.2cm]
&\le3(n+1)\sum_{k=0}^\be 2^{-k}s_{[2^{-k}n]}.
\end{align*}
Hence, $\s_n\le3s_n(A)$, $n\ge0$, and so 
$$
\Psi_\fI\Big(\{\s_n\}_{n\ge0}\Big)\le3\Psi\big(\{\s_n(A)\}_{n\ge0}\big)=
3\|A\|_\fI
\le3c\,\Psi_\fI\Big(\{s_n\}_{n\ge0}\Big).
$$

Let us prove now that (ii)$\imp$(i). Let $\{s_n\}_{n\ge0}$ be a nonincreasing sequence of nonnegative numbers. Put
$$
\xi_n\df\frac1{n+1}\sum_{k=0}^n\s_n=\frac1{n+1}\sum_{k=0}^n\left(\frac1{k+1}\sum_{j=0}^ks_j\right)
=\frac1{n+1}\sum_{j=0}^n\left(\sum_{k=j}^n\frac1{k+1}\right)s_j.
$$
For an arbitrary positive integer $d$, we have
\begin{align*}
\xi_n&\ge\frac{s_{[n/d]}}{n+1}\sum_{j=0}^{[n/d]}\left(\sum_{k=j}^n\frac1{k+1}\right)
\ge\frac{s_{[n/d]}}{n+1}\big([n/d]+1\big)\left(\sum_{k=[n/d]}^n\frac1{k+1}\right)\\[.2cm]
&\ge\frac{\big([n/d]+1\big)s_{[n/d]}}{n+1}\int_{[n/d]}^n\frac{dx}{x+1}
\ge\frac{s_{[n/d]}}d\log\frac{n+2}{[n/d]+1}\ge\frac{\log d}ds_{[n/d]}.
\end{align*}
This together with inequality \rf{xi} applied twice yields
$$
\Psi_\fI\big(\big\{s_{[n/d]}\big\}_{n\ge0}\big)\le\frac{d}{\log d}\Psi_\fI\big(\{\xi_n\}_{n\ge0}\big)
\le\const\frac{d}{\log d}\Psi_\fI\big(\{s_n\}_{n\ge0}\big)
$$
for $d\ge2$. Thus $\b_{\fI,d}<d$ for sufficiently large $d$, and so $\b_\fI<1$. $\bl$

\medskip

{\bf Remark.} Suppose that $\fI$ is a {\it normed} ideal and let
$\bs{C}_\fI$ be the best possible constant in inequality \rf{xi}.
It is easy to see from the proof of Theorem \ref{Bo} that
\bay
\label{Ce}
\bs{C}_\fI\le3\sum_{k=0}^\be2^{-k}\b_{\fI,2^k}.
\ey

\medskip

Let $\bS_p$, $0<p<\be$, be the Schatten--von Neumann class of operators $T$ on Hilbert space
such that
$$
\|T\|_{\bS_p}\df\left(\sum_{j\ge0}\big(s_j(T)\big)^p\right)^{1/p}.
$$
This is a normed ideal for $p\ge1$. We denote by $\bS_{p,\be}$, $0<p<\be$, the ideal that consists of operators $T$ on Hilbert space such that
$$
\|T\|_{\bS_{p,\be}}\df\left(\sup_{j\ge0}(1+j)\big(s_j(T)\big)^p\right)^{1/p}.
$$
The quasinorm $\|\cdot\|_{p,\be}$ is not a norm, but it is equivalent to a norm if $p>1$.
It is easy to see that 
$$
\b_{\bS_p}=\b_{\bS_{p,\be}}=\frac1p,\quad 0<p<\be.
$$
Thus $\bS_p$ and $\bS_{p,\be}$ satisfy the hypotheses of Theorem \ref{Bo} for $p>1$.

It follows easily from \rf{Ce} that for $p>1$,
$$
\bs{C}_{\bS_p}\le3\big(1-2^{1/p-1}\big)^{-1}.
$$

Suppose now that $\fI$ is a quasinormed ideal of operators on Hilbert space.
With a nonnegative integer $l$ we associate the ideal $^{(l)}\fI$ 
that consists of all bounded linear operators on Hilbert space and
is equipped with the norm
$$
\Psi_{^{(l)}\fI}(s_0,s_1,s_2,\cdots)=\Psi(s_0,s_1,\cdots,s_l,0,0,\cdots).
$$
It is easy to see that for every bounded operator $T$,
\begin{align*}
\|T\|_{^{(l)}\fI}&=\sup\big\{\|RT\|_\fI:~\|R\|\le1,~\rank R\le l+1\big\}
\\[.2cm]
&=\sup\big\{\|TR\|_\fI:~\|R\|\le1,~\rank R\le l+1\big\}.
\end{align*}

The following fact is obvious.

\begin{lem}
\label{ob}
Let $\fI$ be a quasinormed ideal. Then for all $l\ge0$,
$$
\bs{C}_{^{(l)}\fI}\le\bs{C}_\fI.
$$
\end{lem}

We refer the reader to \cite{GK} and \cite{BS0} for further information on singular values and normed ideals of operators on Hilbert space.

\

\section{\bf Multiple operator integrals}
\setcounter{equation}{0}
\label{koi}

\

{\bf 3.1. Double operator integrals.}
In this subsection we review some aspects of the theory of double operator integrals. Double operator integrals appeared in the paper \cite{DK} by Daletskii and S.G. Krein. In that paper the authors obtained the following formula
$$
\frac{d}{dt}\big(f(A+tK)-f(A)\big)\Big|_{t=0}=
\iint\frac{f(x)-f(y)}{x-y}\,dE_A(x)K\,dE_A(y)
$$  
for a function $f$ of class $C^2(\R)$, and bounded self-adjoint operators $A$ and $K$ ($E_A$ stands for the spectral measure of $A$). However, the beautiful theory of double operator integrals was developed later by Birman and Solomyak in \cite{BS1}, \cite{BS2}, and \cite{BS3}, see also their survey \cite{BS5}.

Let $(\X,E_1)$ and $(\Y,E_2)$ be spaces with spectral measures $E_1$ and $E_2$
on Hilbert spaces $\h_1$ and $\h_2$. Let us first define double operator integrals
\bay
\label{doi}
\int\limits_{\X}\int\limits_{\Y}\Phi(x,y)\,d E_1(x)\,Q\,dE_2(y),
\ey
for bounded measurable functions $\Phi$ and operators $Q:\h_2\to\h_1$
of Hilbert--Schmidt class $\bS_2$. Consider the set function $F$ whose values are orthogonal
projections on the Hilbert space $\bS_2(\h_2,\h_1)$ of Hilbert--Schmidt operators from $\h_2$ to $\h_1$, which is defined 
on measurable rectangles by
$$
F(\D_1\times\D_2)Q=E_1(\D_1)QE_2(\D_2),\quad Q\in\bS_2(\h_2,\h_1),
$$ 
$\D_1$ and $\D_2$ being measurable subsets of $\X$ and $\Y$. Note that left multiplication by $E_1(\D_1)$
obviously commutes with right multiplication by $E_2(\D_2)$.

It was shown in \cite{BS4} that $F$ extends to a spectral measure on 
$\X\times\Y$. If $\Phi$ is a bounded measurable function on $\X\times\Y$, we define
$$
\int\limits_{\X}\int\limits_{\Y}\Phi(x,y)\,d E_1(x)\,Q\,dE_2(y)=
\left(\,\,\int\limits_{\X_1\times\X_2}\Phi\,dF\right)Q.
$$
Clearly,
$$
\left\|\,\,\int\limits_{\X}\int\limits_{\Y}\Phi(x,y)\,dE_1(x)\,Q\,dE_2(y)\right\|_{\bS_2}
\le\|\Phi\|_{L^\be}\|Q\|_{\bS_2}.
$$

If the transformer
$$
Q\mapsto\int\limits_{\X}\int\limits_{\Y}\Phi(x,y)\,d E_1(x)\,Q\,dE_2(y)
$$
maps the trace class $\bS_1$ into itself, we say that $\Phi$ is a {\it Schur multiplier of $\bS_1$ associated with 
the spectral measures $E_1$ and $E_2$}. In
this case the transformer
\bay
\label{tra}
Q\mapsto\int\limits_{\Y}\int\limits_{\X}\Phi(x,y)\,d E_2(y)\,Q\,dE_1(x),\quad Q\in \bS_2(\h_1,\h_2),
\ey
extends by duality to a bounded linear transformer on the space of bounded linear operators from $\h_1$ to $\h_2$
and we say that the function $\Psi$ on $\X_2\times\X_1$ defined by 
$$
\Psi(y,x)=\Phi(x,y)
$$
is {\it a Schur multiplier of the space of bounded linear operators associated with $E_2$ and $E_1$}.
We denote the space of such Schur multipliers by $\fM(E_2,E_1)$. We also use the notation $\fM(E)\df\fM(E,E)$.

To state a very important formula by Birman and Solomyak, we consider for a continuously differential function
$f$ on $\R$, the divided difference $\dg f$,
$$
(\dg f)(x,y)\df\frac{f(x)-f(y)}{x-y},\quad x\ne y,\quad (\dg f)(x,x)\df f'(x)\quad x,\,y\in\R.
$$
Birman in Solomyak proved in \cite{BS3} that if
 $A$ is a self-adjoint operator (not necessarily bounded),
$K$ is a bounded self-adjoint operator, and
$f$ is a continuously differentiable 
function on $\R$ such that
$\dg f\in\fM(E_{A+K},E_A)$, then
\bay
\label{BSF}
f(A+K)-f(A)=\iint\limits_{\R\times\R}\big(\dg f\big)(x,y)\,dE_{A+K}(x)K\,dE_A(y)
\ey
and
$$
\|f(A+K)-f(A)\|\le\const\|\dg f\|_{\fM}\|K\|,
$$
where $\|\dg f\|_{\fM}$ is the norm of $\dg f$ in $\fM(E_{A+K},E_A)$. 
Here we use the notation $E_A$ for the spectral measure of $A$.


A similar formula and
similar results also hold for unitary operators, in which case we have to integrate the divided
difference $\dg f$ of a function $f$ on the unit circle with respect to the spectral measures of the corresponding operator integrals.

It is easy to see that if a function $\Phi$ on $\X\times\Y$ belongs to the {\it projective tensor
product}
$L^\be(E_1)\hat\otimes L^\be(E_2)$ of $L^\be(E_1)$ and $L^\be(E_2)$ (i.e., $\Phi$ admits a representation
\bay
\label{ptp}
\Phi(x,y)=\sum_{n\ge0}\f_n(x)\psi_n(y),
\ey
where $\f_n\in L^\be(E_1)$, $\psi_n\in L^\be(E_2)$, and
\bay
\label{ptpn}
\sum_{n\ge0}\|\f_n\|_{L^\be}\|\psi_n\|_{L^\be}<\be),
\ey
then $\Phi\in\fM(E_1,E_2)$, i.e., $\Phi$ is a Schur multiplier of the space of bounded linear operators. For such functions $\Phi$ we have
$$
\int\limits_\X\int\limits_\Y\Phi(x,y)\,d E_1(x)Q\,dE_2(y)=
\sum_{n\ge0}\left(\,\int\limits_\X \f_n\,dE_1\right)Q\left(\,\int\limits_\Y \psi_n\,dE_2\right).
$$ 
Note that if $\Phi$ belongs to the projective tensor product $L^\be(E_1)\hat\otimes L^\be(E_2)$, its norm in $L^\be(E_1)\hat\otimes L^\be(E_2)$
is, by definition, the infimum of the  left-hand side of \rf{ptpn} over all representations \rf{ptp}.


More generally, $\Phi$ is a Schur multiplier  if $\Phi$ 
belongs to the {\it integral projective tensor product} $L^\be(E_1)\hat\otimes_{\rm i}L^\be(E_2)$ of $L^\be(E_1)$ and $L^\be(E_2)$, i.e., $\Phi$ admits a representation
$$
\Phi(x,y)=\int_\O \f(x,\o)\psi(y,\o)\,d\s(\o),
$$
where $(\O,\s)$ is a measure space, $\f$ is a measurable function on $\X\times \O$,
$\psi$ is a measurable function on $\Y\times \O$, and
$$
\int_\O\|\f(\cdot,\o)\|_{L^\be(E_1)}\|\psi(\cdot,\o)\|_{L^\be(E_2)}\,d\s(\o)<\be.
$$
If $\Phi\in L^\be(E_1)\hat\otimes_{\rm i}L^\be(E_2)$, then
\bay
\label{iptr}
\iint\limits_{\X\times\Y}\!\Phi(x,y)\,d E_1(x)\,Q\,dE_2(y)\!=\!\!
\int\limits_\O\!\left(\,\int\limits_\X \f(x,\o)\,dE_1(x)\!\right)\!Q\!
\left(\,\int\limits_\Y \psi(y,\o)\,dE_2(y)\!\right)\!d\s(\o).
\ey
Clearly, the function 
$$
\o\mapsto \left(\,\int\limits_\X \f(x,\o)\,dE_1(x)\right)Q
\left(\,\int\limits_\Y \psi(y,\o)\,dE_2(y)\right)
$$
is weakly measurable and
$$
\int\limits_\O\left\|\left(\,\int\limits_\X \f(x,\o)\,dE_1(x)\right)Q
\left(\,\int\limits_\Y \psi(y,\o)\,dE_2(y)\right)\right\|\,d\s(\o)<\be.
$$

It turns out that all Schur multipliers of the space of bounded linear operators can be obtained in this way (see \cite{Pe1}).

This together with the Birman--Solomyak formula \rf{BSF}
implies that if $A$ is a self-adjoint operator and $K$ is a self-adjoint operator that belong to a normed ideal $\fI$, then $f(A+K)-f(A)\in\fI$ and
\bay
\label{BSi}
\|f(A+K)-f(A)\|_\fI
\le\const\|\dg f\|_{L^\be(E_{A+K})\hat\otimes L^\be(E_A)}\|K\|_\fI.
\ey

In connection with the Birman--Solomyak formula, it is important to obtain sharp estimates of divided differences in integral projective tensor products of $L^\be$ spaces. It was shown in \cite{Pe2} that if $f$ is a trigonometric polynomial of degree $d$, then 
\bay
\label{Bp}
\big\|\dg f\big\|_{C(\T)\hat\otimes C(\T)}\le\const d\,\|f\|_{L^\be}.
\ey
On the other hand, it was shown in \cite{Pe4} that if $f$ is a bounded function on $\R$ whose Fourier transform is supported on $[-\s,\s]$
(in other words, $f$ is an entire function of exponential type at most $\s$ that is bounded on $\R$), then $\dg f\in L^\be\hat\otimes_{\rm i}L^\be$
and
\bay
\label{Be}
\big\|\dg f\big\|_{L^\be\hat\otimes_{\rm i} L^\be}\le\const \s\|f\|_{L^\be(\R)}.
\ey
Note that inequalities \rf{Bp} and \rf{Be} were proved in \cite{Pe2} and \cite{Pe4} under the assumption that the Fourier transform of $f$ is supported
on $\Z_+$ (or $\R_+$); however it is very easy to deduce the general results from those partial cases.


\medskip

{\bf 3.2. Multiple operator integrals.}  The approach by Birman and Solomyak to double operator integrals does not generalize to the case of
multiple operator integrals. However, formula \rf{iptr} suggests an approach to multiple operator integrals that is based on integral projective tensor products. This approach was given in \cite{Pe7}.

To simplify the notation, we consider here the case of triple operator integrals; the case of arbitrary multiple operator integrals can be treated in the same way.

Let $(\X,E_1)$, $(\Y,E_2)$, and $(\cZ,E_3)$
be spaces with spectral measures $E_1$, $E_2$, and $E_3$ on Hilbert spaces $\h_1$, $\h_2$, and $\h_3$. Suppose that
$\Phi$ belongs to the integral projective tensor product
$L^\be(E_1)\hat\otimes_{\rm i}L^\be(E_2)\hat\otimes_{\rm i}L^\be(E_3)$, i.e., $\Phi$ admits a representation
\bay
\label{ttp}
\Phi(x,y,z)=\int_\O \f(x,\o)\psi(y,\o)\chi(z,\o)\,d\s(\o),
\ey
where $(\O,\s)$ is a measure space, $\f$ is a measurable function on $\X\times \O$,
$\psi$ is a measurable function on $\Y\times \O$, $\chi$ is a measurable function on $\cZ\times \O$,
and
$$
\int_\O\|\f(\cdot,\o)\|_{L^\be(E)}\|\psi(\cdot,\o)\|_{L^\be(F)}\|\chi(\cdot,\o)\|_{L^\be(G)}\,d\s(\o)<\be.
$$

Suppose now that $T_1$ is a bounded linear operator from $\h_2$ to $\h_1$ and $T_2$ is a bounded linear operator from $\h_3$ to $\h_2$. For a function $\Phi$ in
$L^\be(E_1)\hat\otimes_{\rm i}L^\be(E_2)\hat\otimes_{\rm i}L^\be(E_3)$ of the form \rf{ttp}, we put
\begin{align}
\label{opr}
&\int\limits_\X\int\limits_\Y\int\limits_\cZ\Phi(x,y,z)
\,d E_1(x)T_1\,dE_2(y)T_2\,dE_3(z)\\[.2cm]
\df&\int\limits_\O\left(\,\int\limits_\X \f(x,\o)\,dE_1(x)\right)T_1
\left(\,\int\limits_\Y \psi(y,\o)\,dE_2(y)\right)T_2
\left(\,\int\limits_\cZ \chi(z,\o)\,dE_3(z)\right)\,d\s(\o).\nonumber
\end{align}

It was shown in \cite{Pe7} (see also \cite{ACDS} for a different proof)  that the above definition does not depend on the choice of a representation \rf{ttp}.

It is easy to see that the following inequality holds
$$
\left\|\int\limits_\X\int\limits_\Y\int\limits_\cZ\Phi(x,y,z)
\,dE_1(x)T_1\,dE_2(y)T_2\,dE_3(z)\right\|
\le\|\Phi\|_{L^\be\hat\otimes_{\rm i}L^\be\hat\otimes_{\rm i}L^\be}\|T_1\|\cdot\|T_2\|.
$$

In particular, the triple operator integral on the left-hand side of \rf{opr} can be defined if $\Phi$ belongs to the projective
tensor product $L^\be(E_1)\hat\otimes L^\be(E_2)\hat\otimes L^\be(E_3)$, i.e., $\Phi$ admits a representation
$$
\Phi(x,y,z)=\sum_{n\ge1}\f_n(x)\psi_n(y)\chi_n(z),
$$
where $\f_n\in L^\be(E_1)$, $\psi_n\in L^\be(E_2)$, $\chi_n\in L^\be(E_3)$ and
$$
\sum_{n\ge1}\|\f_n\|_{L^\be(E_1)}\|\psi_n\|_{L^\be(E_2)}\|\chi_n\|_{L^\be(E_3)}<\be.
$$
It is easy to see that if $T_1\in\bS_p$ and $T_2\in\bS_q$, and $1/p+1/q\le1$,
then the triple operator integral \rf{opr} belongs to $\bS_r$ and
$$
\left\|\int\limits_\X\int\limits_\Y\int\limits_\cZ\Phi(x,y,z)
\,dE_1(x)T_1\,dE_2(y)T_2\,dE_3(z)\right\|_{\bS_r}
\le\|\Phi\|_{L^\be\hat\otimes_{\rm i}L^\be\hat\otimes_{\rm i}L^\be}
\|T_1\|_{\bS_p}\cdot\|T_2\|_{\bS_q},
$$
where $1/r=1/p+1/q$.

In a similar way one can define multiple operator integrals, see \cite{Pe7}.

Recall that multiple operator integrals were considered earlier in \cite{Pa} and \cite{S}. However, in those papers the class of functions 
$\Phi$ for which the left-hand side of \rf{opr} was defined is much narrower than in the definition given above.

Multiple operator integrals are used in \cite{Pe7} in connection with the problem of evaluating higher order operator derivatives. 
To obtain formulae for higher order operator derivatives, one has to integrate divided differences of higher orders (see \cite{Pe7}). 

In this paper we are going to integrate divided differences of higher orders to estimate the norms of higher order operator differences \rf{hod}.

For a function $f$ on the circle the  {\it divided differences $\dg^k f$ of order $k$} are defined inductively as follows:
$$
\dg^0f\df f;
$$
if $k\ge1$, then in the case when $\l_1,\l_2,\cdots,\l_{k+1}$ are distinct points in $\T$,
$$
(\dg^{k}f)(\l_1,\cdots,\l_{k+1})\df
\frac{(\dg^{k-1}f)(\l_1,\cdots,\l_{k-1},\l_k)-
(\dg^{k-1}f)(\l_1,\cdots,\l_{k-1},\l_{k+1})}{\l_{k}-\l_{k+1}}
$$
(the definition does not depend on the order of the variables). Clearly,
$$
\dg f=\dg^1f.
$$
If $f\in C^k(\T)$, then $\dg^{k}f$ extends by continuity to a function defined for all points $\l_1,\l_2,\cdots,\l_{k+1}$.

It can be shown that
$$
({\frak D}^n\var)(\l_1,\dots,\l_{n+1})=\sum\limits_{k=1}^{n+1}\var(\l_k)
\prod\limits_{j=1}^{k-1}(\l_k-\l_j)^{-1}\prod\limits_{j=k+1}^{n+1}(\l_k-\l_j)^{-1}.
$$

Similarly, one can define the divided difference of order $k$ for functions on the real line.

It was shown in \cite{Pe7} that if $f$ is a trigonometric polynomial of degree $d$, then
\bay
\label{Bok}
\big\|\dg^k f\big\|_{C(\T)\hat\otimes\cdots\hat\otimes C(\T)}\le\const d^k\|f\|_{L^\be}.
\ey
It was also shown in \cite{Pe7} that if $f$ is an entire function of exponential type at most $\s$ and is bounded on $\R$, then
\bay
\label{Boke}
\big\|\dg^k f\big\|_{L^\be\hat\otimes_{\rm i}\cdots\hat\otimes_{\rm i} L^\be}\le\const \s^k\|f\|_{L^\be(\R)}.
\ey


\medskip

{\bf 3.3. Multiple operator integrals with respect to semi-spectral measures.}
Let $\h$ be a Hilbert space and let $(\X,{\frak B})$ be a measurable space.
A map $\E$ from ${\frak B}$ to the algebra $\B(\h)$ of all bounded operators on $\h$ is called a {\it semi-spectral measure}
if 
$$
\E(\D)\ge\0,\quad\D\in{\frak B},
$$
$$
\E(\varnothing)=\0\quad\mbox{and}\quad\E(\X)=I,
$$
and for a sequence $\{\D_j\}_{j\ge1}$ of disjoint sets in ${\frak B}$,
$$
\E\left(\bigcup_{j=1}^\be\D_j\right)=\lim_{N\to\be}\sum_{j=1}^N\E(\D_j)\quad\mbox{in the weak operator topology}.
$$

\medskip

If $\K$ is a Hilbert space, $(\X,{\frak B})$ is a measurable space,  $E:{\frak B}\to \B(\K)$ is a spectral measure, and $\h$ is
a subspace of $\K$, then it is easy to see that the map $\E:{\frak B}\to \B(\h)$ defined by
\bay
\label{dil}
\E(\D)=P_\h E(\D)\big|\h,\quad\D\in{\frak B},
\ey
is a semi-spectral measure. Here $P_\h$ stands for the orthogonal projection onto $\h$.

Naimark proved in \cite{Na}  that all semi-spectral measures can be obtained in this way, i.e.,
a semi-spectral measure is always a {\it compression} of a spectral measure. A spectral measure $E$ satisfying \rf{dil} is called a {\it spectral dilation of the semi-spectral measure} $\E$.

A spectral dilation $E$ of a semi-spectral measure $\E$ is called {\it minimal} if 
$$
\K=\clos\spn\{E(\D)\h:~\D\in{\frak B}\}.
$$

It was shown in \cite{MM} that if $E$ is a minimal spectral dilation of a semi-spectral measure $\E$, then
$E$ and $\E$ are mutually absolutely continuous and all minimal spectral dilations of a semi-spectral measure are isomorphic in the natural sense.

If $\f$ is a bounded complex-valued measurable function on $\X$ and $\E:{\frak B}\to \B(\h)$ is a semi-spectral measure, then the integral
\bay
\label{iss}
\int_\X \f(x)\,d\E(x)
\ey
can be defined as
\bay
\label{voi}
\int_\X \f(x)\,d\E(x)=\left.P_\h\left(\int_\X \f(x)\,d E(x)\right)\right|\h,
\ey
where $E$ is a spectral dilation of $\E$. It is easy to see that the right-hand side of \rf{voi} does not depend on the choice
of a spectral dilation. The integral \rf{iss} can also be computed as the limit of sums
$$
\sum \f(x_\a)\E(\D_\a),\quad x_\a\in\D_\a,
$$
over all finite measurable partitions $\{\D_\a\}_\a$ of $\X$.

If $T$ is a contraction on a Hilbert space $\h$, then by the Sz.-Nagy dilation theorem
(see \cite{SNF}),  $T$ has a unitary dilation, i.e., there exist a Hilbert space $\K$ such that
$\h\subset\K$ and a unitary operator $U$ on $\K$ such that
\bay
\label{DT}
T^n=P_\h U^n\big|\h,\quad n\ge0,
\ey
where $P_\h$ is the orthogonal projection onto $\h$. Let $E_U$ be the spectral measure of $U$.
Consider the operator set function $\E$ defined on the Borel subsets of the unit circle $\T$ by
$$
\E(\D)=P_\h E_U(\D)\big|\h,\quad\D\subset\T.
$$
Then $\E$ is a semi-spectral measure. It follows immediately from
\rf{DT} that 
\bay
\label{step}
T^n=\int_\T \z^n\,d\E(\z)=P_\h\int_\T\z^n\,dE_U(\z)\Big|\h,\quad n\ge0.
\ey
Such a semi-spectral measure $\E$ is called a {\it semi-spectral measure} of $\T$.
Note that it is not unique. To have uniqueness, we can consider a minimal unitary dilation $U$ of $T$,
which is unique up to an isomorphism (see \cite{SNF}).

It follows easily from  \rf{step} that 
$$
f(T)=P_\h\int_\T f(\z)\,dE_U(\z)\Big|\h
$$
for an arbitrary function $\f$ in the disk-algebra $C_A$.

In \cite{Pe3} and  \cite{Pe8} double operator integrals and multiple operator integrals with respect to semi-spectral measures were introduced.

Suppose that $(\X_1,{\frak B}_1)$ and $(\X_2,{\frak B}_2)$ are measurable spaces, and
$\E_1:{\frak B}_1\to \B(\h_1)$ and $\E_2:{\frak B}_2\to \B(\h_2)$ are semi-spectral measures.
Then double operator integrals
$$
\iint\limits_{\X_1\times\X_2}\Phi(x_1,x_2)\,d\E_1(x_1)Q\,d\E_2(X_2).
$$
were defined in \cite{Pe8} in the case when $Q\in\bS_2$ and $\Phi$ is a bounded measurable function. Double operator integrals were also defined in \cite{Pe8} in the case when $Q$ is a bounded linear operator and $\Phi$ belongs to the integral projective tensor product of the spaces $L^\be(\E_1)$
and $L^\be(\E_2)$.

In particular, the following analog of the Birman--Solomyak formula holds:
\bay
\label{BSc}
f(R)-f(T)=\iint\limits_{\T\times\T}\big(\dg f\big)(\z,\t)\,d\E_R(\z)(R-T)\,d\E_T(\t).
\ey
Here $T$ and $R$ contractions on Hilbert space, $\E_T$ and $\E_R$ are their semi-spectral measures, and $f$ is an analytic  function in $\dd$ of class
$\big(B_{\be1}^1\big)_+$.

Similarly, multiple operator integrals with respect to semi-spectral measures were defined in \cite{Pe8} for functions that belong to the integral projective tensor product of the corresponding $L^\be$ spaces.

\

\section{\bf Self-adjoint operators. Sufficient conditions}
\setcounter{equation}{0}
\label{sa}

\

For $l\ge0$ and $p>0$, we consider the normed ideal $\bS_p^l$ that consists of all bounded linear operators equipped with the norm
$$
\|T\|_{\bS_p^l}\df\left(\sum_{j=0}^l\big(s_j(T)\big)^p\right)^{1/p}.
$$
It is well known that $\|\cdot\|_{\bS_p^n}$ is a norm for $p\ge1$ (see \cite{BS0}).
Note that $\bS_p^l={^{(l)}\bS_p}$, see \S\,\ref{ide}. 

\begin{thm}
\label{gl}
Let $0<\a<1$. Then there exists a positive number $c>0$ such that for every
$l\ge0$,  $p\in[1,\be)$,  $f\in\L_\a(\R)$, and for arbitrary self-adjoint operators $A$ and $B$ on Hilbert space with bounded $A-B$, the following inequality holds:
$$
s_j\big(f(A)-f(B)\big)\le c\,\|f\|_{\L_\a(\R)}(1+j)^{-\a/p}\|A-B\|_{\bS_p^l}^\a
$$
for every $j\le l$.
\end{thm}

\Pf Put $f_n\df f*W_n+f*W_n^\sharp$, $n\in\Z$, and fix an integer $N$. We have
by \rf{BSi} and \rf{Be},
\begin{align*}
\left\|\sum_{n=-\be}^N\big(f_n(A)-f_n(B)\big)\right\|_{\bS_p^l}
&\le\sum_{n=-\be}^N\big\|f_n(A)-f_n(B)\big\|_{\bS_p^l}\\[.2cm]
&\le\const\sum_{n=-\be}^N2^n\|f_n\|_{L^\be}\|A-B\|_{\bS_p^l}\\[.2cm]
&\le\const\|f\|_{\L_\a(\R)}\sum_{n=-\be}^N2^{n(1-\a)}\|A-B\|_{\bS_p^l}\\[.2cm]
&\le\const2^{N(1-\a)}\|f\|_{\L_\a(\R)}\|A-B\|_{\bS_p^l}.
\end{align*}
On the other hand,
\begin{align*}
\left\|\sum_{n>N}\big(f_n(A)-f_n(B)\big)\right\|
&\le2\sum_{n>N}\|f_n\|_{L^\be}\\[.2cm]
&\le\const\|f\|_{\L_\a(\R)}\sum_{n>N}2^{-n\a}
\le\const2^{-N\a}\|f\|_{\L_\a(\R)}.
\end{align*}
Put
$$
R_N\df\sum_{n=-\be}^N\big(f_n(A)-f_n(B)\big)
\qm Q_N\df\sum_{n>N}\big(f_n(A)-f_n(B)\big).
$$
Clearly, for $j\le l$,
\begin{align*}
s_j\big(f(A)-f(B)\big)&\le s_j(R_N)+\|Q_N\|
\le(1+j)^{-1/p}\|f(A)-f(B)\|_{\bS_p^l}+\|Q_N\|\\[.2cm]
&\le\const\left((1+j)^{-1/p}2^{N(1-\a)}\|f\|_{\L_\a(\R)}\|A-B\|_{\bS_p^l}
+2^{-N\a}\|f\|_{\L_\a(\R)}\right).
\end{align*}
To obtain the desired estimate, it suffices to choose the number $N$ so that
$$
2^{-N}<(1+j)^{-1/p}\|A-B\|_{\bS_p^l}\le2^{-N+1}.\quad\bl
$$

\begin{thm}
\label{S1}
Let $0<\a<1$. Then there exists a positive number $c>0$ such that for every
$f\in\L_\a(\R)$ and arbitrary self-adjoint operators $A$ and $B$ on Hilbert space with $A-B\in\bS_1$, the operator $f(A)-f(B)$ belongs to 
$\bS_{\frac1\a,\be}$ and the following inequality holds:
$$
\big\|f(A)-f(B)\big\|_{\bS_{\frac1\a,\be}}\le c\,\|f\|_{\L_\a(\R)}
\|A-B\|_{\bS_1}^\a.
$$
\end{thm}

\Pf This is an immediate consequence of Theorem \ref{gl} in the case $p=1$. $\bl$

Note that the assumptions of Theorem \ref{S1} do not imply that 
$f(A)-f(B)\in\bS_{1/\a}$. In \S\,\ref{frsa} we obtain a necessary condition on 
$f$ for $f(A)-f(B)\in\bS_{1/\a}$ whenever $A-B\in\bS_1$. 

The following result ensures that the assumption that $A-B\in\bS_1$
implies that \lb$f(A)-f(B)\in\bS_{1/\a}$ under a slightly more restrictive condition on $f$.

\begin{thm}
\label{S1s}
Let $0<\a\le1$. Then there exists a positive number $c>0$ such that for every
$f\in B_{\be1}^\a(\R)$ and arbitrary self-adjoint operators $A$ and $B$ on Hilbert space with $A-B\in\bS_1$, the operator $f(A)-f(B)$ belongs to 
$\bS_{1/\a}$ and the following inequality holds:
$$
\big\|f(A)-f(B)\big\|_{\bS_{1/\a}}\le c\,\|f\|_{B_{\be1}^\a(\R)}
\|A-B\|_{\bS_1}^\a.
$$
\end{thm}

Note that in the case $\a=1$ this was proved earlier in \cite{Pe4}.

{\bf Proof of Theorem \ref{S1s}.}
 Put $f_n=f*W_n+f*W_n^\sharp$. Clearly, $f_n$ is trace class perturbations preserving
and it is easy to see that
\bay
\label{Had}
\|f_n(A)-f_n(B)\|_{\bS_{1/\a}}\le
\|f_n(A)-f_n(B)\|_{\bS_1}^\a\|f_n(A)-f_n(B)\|^{1-\a}.
\ey
Since $f(A)-f(B)=\sum\limits_{n\in\Z}\big(f_n(A)-f_n(B)\big)$, it suffices to prove that
$$
\sum_{n\in\Z}\big\|f_n(A)-f_n(B)\big\|_{\bS_{1/\a}}<\be.
$$
We have by \rf{Had} and \rf{BSi},
\begin{align*}
\sum_{n\in\Z}\big\|f_n(A)-f_n(B)\big\|_{\bS_{1/\a}}&\le
\sum_{n\in\Z}\big\|f_n(A)-f_n(B)\big\|_{\bS_1}^\a\cdot
\big\|f_n(A)-f_n(B)\big\|^{1-\a}\\[.2cm]
&\le\const\sum_{n\in\Z}2^{n\a}\|f_n\|_{L^\be}^\a
\cdot2^{1-\a}\|f_n\|_{L^\be}^{1-\a}\|A-B\|_{\bS_1}^\a\\[.2cm]
&\le\const\sum_{n\in\Z}2^{n\a}\|f_n\|_{L^\be}\|A-B\|_{\bS_1}^\a\\[.2cm]
&\le\const\|f\|_{B_{\be1}^\a(\R)}\|A-B\|_{\bS_1}^\a.
\quad\bl
\end{align*}

\begin{thm}
\label{sigma}
Let $0<\a<1$. Then there exists a positive number $c>0$ such that for every
$f\in\L_\a(\R)$ and arbitrary self-adjoint operators $A$ and $B$ on Hilbert space with bounded $A-B$, the following inequality holds:
$$
s_j\Big(\big|f(A)-f(B)\big|^{1/\a}\Big)\le c\,\|f\|_{\L_\a(\R)}^{1/\a}
\s_j(A-B),\quad j\ge0.
$$
\end{thm}

\Pf It suffices to apply Theorem \ref{gl} with $l=j$ and $p=1$. $\bl$

Now we are in a position to obtain a general result in the case 
$f\in\L_\a(\R)$ and $A-B\in\fI$ for an arbitrary quasinormed ideal $\fI$ with upper Boyd index less than 1.

\begin{thm}
\label{osn}
Let $0<\a<1$. Then there exists a positive number $c>0$ such that for every
$f\in\L_\a(\R)$, for an arbitrary quasinormed ideal $\fI$ with $\b_\fI<1$, and for arbitrary self-adjoint operators $A$ and $B$ on Hilbert space with $A-B\in\fI$,
the operator $\big|f(A)-f(B)\big|^{1/\a}$ belongs to $\fI$
and the following inequality holds:
$$
\Big\|\big|f(A)-f(B)\big|^{1/\a}\Big\|_\fI\le c\,\bs{C}_\fI
\|f\|_{\L_\a(\R)}^{1/\a}\|A-B\|_\fI.
$$
\end{thm}

\Pf The result follows from Theorems \ref{sigma} and \ref{Bo}. $\bl$

We can reformulate Theorem \ref{osn} in the following way.

\begin{thm}
\label{pere}
Under the hypothesis of Theorem {\em\ref{osn}}, the operator $f(A)-f(B)$
belongs to $\fI^{\{1/\a\}}$ and
$$
\big\|f(A)-f(B)\big\|_{\fI^{\{1/\a\}}}\le c^\a\,\bs{C}_\fI^\a\|f\|_{\L_\a(\R)}
\|A-B\|^\a_\fI.
$$
\end{thm}

We deduce now some more consequences of Theorem \ref{osn}.

\begin{thm}
\label{spl}
Let $0<\a<1$ and $1<p<\be$. Then there exists a positive number $c$ such that
for every $f\in\L_\a(\R)$, every $l\in\Z_+$, and arbitrary self-adjoint operators $A$ and $B$ with bounded $A-B$, the following inequality holds:
$$
\sum_{j=0}^l\left(s_j\Big(\big|f(A)-f(B)\big|^{1/\a}\Big)\right)^p\le
c\,\|f\|_{\L_\a(\R)}^{p/\a}\sum_{j=0}^l\big(s_j(A-B)\big)^p.
$$
\end{thm}

\Pf The result immediately follows from Theorem \ref{osn} and Lemma \ref{ob}.
$\bl$

\begin{thm}
\label{sp}
Let $0<\a<1$ and $1<p<\be$. Then there exists a positive number $c$ such that
for every $f\in\L_\a(\R)$ and for arbitrary self-adjoint operators $A$ and $B$ with $A-B\in\bS_p$, the operator $f(A)-f(B)$ belongs to $\bS_{p/\a}$ and
the following inequality holds:
$$
\big\|f(A)-f(B)\big\|_{\bS_{p/\a}}\le c\,\|f\|_{\L_\a(\R)}\|A-B\|^\a_{\bS_p}.
$$
\end{thm}

\Pf The result is an immediate consequence of Theorem \ref{spl}. $\bl$

To proceed to higher order differences, we need the following well-known inequality:
\bay
\label{rpq}
\|T_1T_2\|_{\bS_r^l}\le\|T_1\|_{\bS_p^l}\|T_2\|_{\bS_q^l},
\ey
where $T_1$ and $T_2$ bounded operator on Hilbert space and $1/p+1/q\le1$.
Inequality \rf{rpq} can be deduced from the corresponding inequality for
$\bS_p$ norms. Indeed, let $R$ be an operator of rank $l$ such that 
$\|T_1T_2\|_{\bS_r^l}=\|T_1T_2R\|_{\bS_r}$. There exists an orthogonal projection $P$ of rank $l$ such that $\|T_1T_2R\|_{\bS_r}=\|PT_1T_2R\|_{\bS_r}$. Then
$$
\|T_1T_2\|_{\bS_r^l}=\|PT_1T_2R\|_{\bS_r}\le\|PT_1\|_{\bS_p}\|T_2R\|_{\bS_q}\le
\|T_1\|_{\bS_p^l}\|T_2\|_{\bS_q^l}.
$$

Suppose now that $m-1\le\a<m$ and $f\in\L_\a(\R)$. For a self-adjoint operator $A$ and a bounded self-adjoint operator $K$, we consider the finite difference
$$
\big(\D_K^mf\big)(A)\df\sum_{j=0}^m(-1)^{m-j}\left(\begin{matrix}m\\j\end{matrix}\right)f\big(A+jK\big).
$$
In the case when $A$ is unbounded, by the right-hand side we mean the following operator
$$
\sum_{n\in\Z}\sum_{j=0}^m(-1)^{m-j}
\left(\begin{matrix}m\\j\end{matrix}\right)f_n\big(A+jK\big),
$$
where as usual, $f_n=f*W_n+f*W_n^\sharp$. It has been proved in \cite{AP2} that
under the above assumptions, 
$$
\sum_{n\in\Z}\left\|\sum_{j=0}^m(-1)^{m-j}
\left(\begin{matrix}m\\j\end{matrix}\right)f_n\big(A+jK\big)\right\|<\be.
$$
(We refer the reader to \cite{AP3}, where the situation with unbounded $A$ will be discussed in detail.)

We are going to use the following representation for $\big(\D_K^mf\big)(A)$ in terms of multiple operator integrals:
\begin{align}
\label{moi}
\big(\D_K^m&f\big)(A)=\\[.2cm]
&m!\underbrace{\int\cdots\int}_{m+1}(\dg^{m}f)(x_1,\cdots,x_{m+1})
\,dE_A(x_1)K\,dE_{A+K}(x_2)K\cdots K\,dE_{A+mK}(x_{m+1}),\nonumber
\end{align}
where $A$ is a self-adjoint operator, $K$ is a bounded self-adjoint operator,
and $f\in B_{\be1}^m(\R)$. Formula \rf{moi} was obtained in \cite{AP2}.

It follows from \rf{moi}, \rf{Boke}, and \rf{rpq} that if $p\ge m\ge1$, $l\ge0$, and $f$ is an entire function of exponential type at most $\s$ that is bounded on $\R$, then
\bay
\label{spln}
\big\|\big(\D_K^mf\big)(A)\big\|_{\bS_{\frac{p}{m}}^l}
\le\const\s^m\|f\|_{L^\be}\|K\|^m_{\bS_p^l}.
\ey
Moreover, the constant in \rf{spln} does not depend on $p$.

Inequality \rf{spln} can be generalized. Suppose that $\fI$ is a normed ideal such that $\fI^{\{1/m\}}$ is also a normed ideal. Suppose that $K\in\fI$. Then
$\big(\D_K^mf\big)(A)\in\fI^{\{1/m\}}$ and
\bay
\label{li}
\big\|\big(\D_K^mf\big)(A)\big\|_{\fI^{\{1/m\}}}
\le\const\s^m\|f\|_{L^\be}\|K\|^m_\fI.
\ey

\begin{thm}
\label{vyp}
Let $\a>0$ and $m-1\le\a<m$. There exists a positive number $c$ such that for every $l\ge0$, $p\in[m,\be)$, $f\in\L_\a(\R)$, and for arbitrary self-adjoint operator $A$ and bounded self-adjoint operator $K$, the following inequality holds:
$$
s_j\Big(\big(\D_K^mf\big)(A)\Big)\le c\,\|f\|_{\L_\a(\R)}(1+j)^{-\a/p}\|K\|^\a_{\bS_p^l}
$$
for $j\le l$.
\end{thm}

\Pf As in the proof of Theorem \ref{gl}, we put
$$
R_N\df\sum_{n\le N}\big(\D_K^mf_n\big)(A)\qm 
Q_N\df\sum_{n>N}\big(\D_K^mf_n\big)(A).
$$
It follows \rf{spln} that
\begin{align*}
\|R_N\|_{\bS^l_{p/m}}&\le\sum_{n\le N}
\const2^{mn}\|f_n\|_{L^\be}\|K\|^m_{\bS^l_{p}}\\[.2cm]
&\le\|K\|^m_{\bS^l_{p}}\|f\|_{\L_\a(\R)}\sum_{n\le N}2^{(m-\a)n}
\le2^{(m-\a)N}\|f\|_{\L_\a(\R)}\|K\|^m_{\bS^l_{p}}.
\end{align*}
On the other hand, it is easy to see that
$$
\|Q_N\|\le\const\sum_{n>N}\|f_n\|_{L^\be}
\le\|f\|_{\L_\a(\R)}\sum_{n>N}2^{-n\a}
\le2^{-\a N}\|f\|_{\L_\a(\R)}.
$$

Hence, 
\begin{align*}
s_j\Big(\big(\D_K^mf_n\big)(A)\Big)&\le s_j(R_N)+\|Q_N\|
\le(1+j)^{-m/p}\|R_N\|_{\bS^l_{p/m}}+\|Q_N\|\\[.2cm]
&\le\const\|f\|_{\L_\a(\R)}\left((1+j)^{-m/p}2^{(m-\a)N}\|K\|^m_{\bS^l_{p}}
+2^{-\a N}\right).
\end{align*}
To complete the proof, it suffices to choose $N$ such that
$$
2^{-N}<(1+j)^{-1/p}\|K\|_{\bS^l_{p}}\le2^{-N+1}.\quad\bl
$$

The following result is an immediate consequence from Theorem \ref{vyp}.

\begin{thm}
\label{vsl}
Let $\a>0$ and $m-1\le\a<m$. There exists a positive number $c$ such that for every $f\in\L_\a(\R)$, and for an arbitrary self-adjoint operator $A$ and an arbitrary self-adjoint operator $K$ of class $\bS_m$, the operator $\big(\D_K^mf\big)(A)$ belongs to
$\bS_{\frac{m}{\a},\be}$ and the following inequality holds:
$$
\big\|\big(\D_K^mf\big)(A)\big\|_{\bS_{\frac{m}{\a},\be}}
\le c\,\|f\|_{\L_\a(\R)}\|K\|_{\bS_m}^\a.
$$
\end{thm}

As in the case $0<\a<1$ (see Theorem \ref{S1s}), we are going to improve the conclusion of Theorem \ref{vsl} under a slightly more restrictive assumption on $f$. Note that in the following theorem $\a$ is allowed to be equal to $m$.

\begin{thm}
\label{p=m}
Let $\a>0$ and $m-1\le\a\le m$. There exists a positive number $c$ such that for every $f\in B_{\be1}^\a(\R)$, and for an arbitrary self-adjoint operator $A$  and
an arbitrary self-adjoint operator $K$ of class $\bS_m$, the operator $\big(\D_K^mf\big)(A)$ belongs to
$\bS_{\frac{m}{\a}}$ and the following inequality holds:
$$
\big\|\big(\D_K^mf\big)(A)\big\|_{\bS_{\frac{m}{\a}}}
\le c\,\|f\|_{B_{\be1}^\a(\R)}\|K\|_{\bS_m}^\a.
$$
\end{thm}

\Pf Clearly,
$$
\big\|\big(\D_K^mf_n\big)(A)\big\|_{\bS_{\frac{m}{\a}}}
\le\big\|\big(\D_K^mf_n\big)(A)\big\|_{\bS_1}^{\a/m}
\big\|\big(\D_K^mf_n\big)(A)\big\|^{1-\a/m}.
$$
By \rf{li},
$$
\big\|\big(\D_K^mf_n\big)(A)\big\|_{\bS_1}\le
\const2^{mn}\|f_n\|_{L^\be}\|K\|_{\bS_m}^m.
$$
Thus
\begin{align*}
\sum_{n\in\Z}\big\|\big(\D_K^mf_n\big)(A)\big\|_{\bS_{\frac{m}{\a}}}&
\le\sum_{n\in\Z}\big\|\big(\D_K^mf_n\big)(A)\big\|_{\bS_1}^{\a/m}
\big\|\big(\D_K^mf_n\big)(A)\big\|^{1-\a/m}\\[.2cm]
&\le\const
\sum_{n\in\Z}
2^{\a n}\|f_n\|_{L^\be}^{\a/m}\|K\|_{\bS_m}^\a\|f_n\|_{L^\be}^{1-\a/m}
\\[.2cm]
&\le\const\|K\|_{\bS_m}^\a\sum_{n\in\Z}2^{\a n}\|f_n\|_{L^\be}
\le\const\|f\|_{B_{\be1}^\a(\R)}\|K\|_{\bS_m}^\a.\quad\bl
\end{align*}

Recall that for a bounded linear operator $T$ the numbers, $\s_j(T)$ are defined by \rf{sin}.

\begin{thm}
\label{vsin}
Let $\a>0$ and $m-1\le\a<m$. There exists a positive number $c$ such that for every $f\in\L_\a(\R)$, and for arbitrary self-adjoint operator $A$ and bounded self-adjoint operator $K$, the following inequality holds:
$$
s_j\Big(\big|\big(\D_K^mf\big)(A)\big|^{m/\a}\Big)\le c\,
\|f\|_{\L_\a(\R)}^{m/\a}\s_j\big(|K|^m\big),\quad j\ge0.
$$
\end{thm}

\Pf The result follows immediately from Theorem \ref{vyp} in the case $j=l$ and $p=m$. $\bl$

\begin{thm}
\label{id}
Let $\a>0$ and $m-1\le\a<m$. There exists a positive number $c$ such that for every $f\in\L_\a(\R)$, every quasinormed ideal $\fI$
with $\b_\fI<m^{-1}$, and for arbitrary self-adjoint operator $A$ and bounded self-adjoint operator $K$, the following inequality holds:
$$
\Big\|\big|\big(\D_K^mf\big)(A)\big|^{1/\a}\Big\|_\fI
\le c\,\bs{C}^{1/m}_{\fI^{\{1/m\}}}\|f\|_{\L_\a(\R)}^{1/\a}\|K\|_\fI.
$$
\end{thm}

\Pf Clearly, $|K|^m\in\fI^{\{1/m\}}$ and $\b_{\fI^{\{1/m\}}}=m\b_\fI<1$.
Therefore, by Theorem \ref{vsin}, 
$$
\Big\|\big|\big(\D_K^mf\big)(A)\big|^{m/\a}\Big\|_{\fI^{\{1/m\}}}
\le c\,\bs{C}^{1/m}_{\fI^{\{1/m\}}}\|f\|_{\L_\a(\R)}^{m/\a}\big\|
|K|^m\big\|_{\fI^{\{1/m\}}}
$$
which implies the result. $\bl$

\begin{thm}
\label{sjm}
Let $\a>0$, $m-1\le\a<m$, and $m<p<\be$. There exists a positive number $c$ such that for every $f\in\L_\a(\R)$, every $l\in\Z_+$, and for arbitrary self-adjoint operator $A$ and bounded self-adjoint operator $K$, the following inequality holds:
$$
\sum_{j=0}^l\left(s_j\left(\big|\big(\D_K^mf\big)(A)\big|^{1/\a}\right)\right)^p
\le c\,\|f\|_{\L_\a(\R)}^{p/\a}\sum_{j=0}^l\big(s_j(K)\big)^p.
$$
\end{thm}

\Pf The result follows from Theorem \ref{id} and Lemma \ref{ob}. $\bl$

The last theorem of this section is an immediate consequence of Theorem \ref{sjm}.

\begin{thm}
\label{posl}
Let $\a>0$, $m-1\le\a<m$, and $m<p<\be$. There exists a positive number $c$ such that for every $f\in\L_\a(\R)$, for an arbitrary self-adjoint operator $A$, and an arbitrary self-adjoint operator $K$ of class $\bS_p$, the following inequality holds:
$$
\left\|\big(\D_K^mf\big)(A)\right\|_{\bS_{p/\a}}
\le c\,\|f\|_{\L_\a(\R)}\|K\|_{\bS_p}^\a.
$$
\end{thm}

\

\section{\bf Unitary operators. Sufficient conditions}
\setcounter{equation}{0}
\label{u}

\

In this section we are going to obtain analogs of the results of the previous section for functions of unitary operators. In the case of first order differences we can use the Birman--Solomyak formula for functions of unitary operators and the proofs are the same as in the case of functions of self-adjoint operators. However, in the case of higher order differences, formulae that express a difference of order $m$ involves not only multiple operator integrals of multiplicity $m+1$, but also multiple operator integrals of lower multiplicities, see \cite{AP2}. This makes proofs more complicated than in the self-adjoint case.

We start with first order differences. If $U$ and $V$ are unitary operators, then
by the Birman--Solomyak formula,
\bay
\label{BSu}
f(U)-f(V)=
\iint\limits_{\T\times\T}\frac{f(\z)-f(\t)}{\z-\t}\,dE_U(\z)(U-V)\,dE_V(\t),
\ey
whenever the divided difference $\dg f$ belongs to $L^\be\hat\otimes L^\be$.
Here $E_U$ and $E_V$ are the spectral measures of $U$ and $V$. Recall that it was shown in \cite{Pe2} that \rf{BSu} holds if $f\in B_{\be1}^1$. 

It follows from \rf{Bp} that if $\fI$ is a normed ideal, $U-V\in\fI$ and $f$ is a trigonometric polynomial of degree $d$, then $f(U)-f(V)\in\fI$ and
\bay
\label{UVI}
\|f(U)-f(V)\|_\fI\le\const d\|f\|_{L^\be}\|U-V\|_\fI.
\ey
Moreover, the constant does not depend on $\fI$.

\begin{thm}
\label{sju}
Let $0<\a<1$. Then there exists a positive number $c>0$ such that for every
$l\ge0$,  $p\in[1,\be)$,  $f\in\L_\a$, and for arbitrary unitary operators $U$ and $V$ on Hilbert space, the following inequality holds:
$$
s_j\big(f(U)-f(V)\big)\le c\,\|f\|_{\L_\a}(1+j)^{-\a/p}\|U-V\|_{\bS_p^l}^\a
$$
for every $j\le l$.
\end{thm}

\begin{thm}
\label{S1u}
Let $0<\a<1$. Then there exists a positive number $c>0$ such that for every
$f\in\L_\a$ and arbitrary unitary operators $U$ and $V$ on Hilbert space with $U-V\in\bS_1$, the operator $f(U)-f(V)$ belongs to 
$\bS_{\frac1\a,\be}$ and the following inequality holds:
$$
\big\|f(U)-f(V)\big\|_{\bS_{\frac1\a,\be}}\le c\,\|f\|_{\L_\a}
\|U-V\|_{\bS_1}^\a.
$$
\end{thm}

As in the self-adjoint case, the assumptions of Theorem \ref{S1u} do not imply that 
\lb$f(U)-f(V)\in\bS_{1/\a}$. In \S\,\ref{fru} we obtain a necessary condition on 
$f$ for $f(U)-f(V)\in\bS_{1/\a}$, whenever $U-V\in\bS_1$.

\begin{thm}
\label{S1Bu}
Let $0<\a\le1$. Then there exists a positive number $c>0$ such that for every
$f\in B_{\be1}^\a$ and arbitrary unitary operators $U$ and $V$ on Hilbert space with $U-V\in\bS_1$, the operator $f(U)-f(V)$ belongs to 
$\bS_{1/\a}$ and the following inequality holds:
$$
\big\|f(U)-f(V)\big\|_{\bS_{1/\a}}\le c\,\|f\|_{B_{\be1}^\a}
\|U-V\|_{\bS_1}^\a.
$$
\end{thm}

Note that in the case $\a=1$ this was proved earlier in \cite{Pe2}.

\begin{thm}
\label{sigu}
Let $0<\a<1$. Then there exists a positive number $c>0$ such that for every
$f\in\L_\a$ and arbitrary unitary operators $U$ and $V$ on Hilbert space, the following inequality holds:
$$
s_j\Big(\big|f(U)-f(V)\big|^{1/\a}\Big)\le c\,\|f\|_{\L_\a}^{1/\a}
\s_j(U-V),\quad j\ge0.
$$
\end{thm}

Recall that the numbers $\s_j(U-V)$ are defined in \rf{sin}.

\begin{thm}
\label{osnu}
Let $0<\a<1$. Then there exists a positive number $c>0$ such that for every
$f\in\L_\a$, for an arbitrary quasinormed ideal $\fI$ with $\b_\fI<1$, and for arbitrary unitary operators $U$ and $V$ on Hilbert space with $U-V\in\fI$,
the operator $\big|f(U)-f(V)\big|^{1/\a}$ belongs to $\fI$
and the following inequality holds:
$$
\Big\|\big|f(U)-f(V)\big|^{1/\a}\Big\|_\fI\le c\,\bs{C}_\fI
\|f\|_{\L_\a(\R)}^{1/\a}\|U-V\|_\fI.
$$
\end{thm}

\begin{thm}
\label{splu}
Let $0<\a<1$ and $1<p<\be$. Then there exists a positive number $c$ such that
for every $f\in\L_\a$, every $l\in\Z_+$, and arbitrary unitary operators $U$ and $V$, the following inequality holds:
$$
\sum_{j=0}^l\left(s_j\Big(\big|f(U)-f(V)\big|^{1/\a}\Big)\right)^p\le
c\,\|f\|_{\L_\a}^{p/\a}\sum_{j=0}^l\big(s_j(U-V)\big)^p.
$$
\end{thm}

\begin{thm}
\label{spu}
Let $0<\a<1$ and $1<p<\be$. Then there exists a positive number $c$ such that
for every $f\in\L_\a$ and for arbitrary unitary operators $U$ and $V$ with $U-V\in\bS_p$, the operator $f(U)-f(V)$ belongs to $\bS_{p/\a}$ and
the following inequality holds:
$$
\big\|f(U)-f(V)\big\|_{\bS_{p/\a}}\le c\,\|f\|_{\L_\a}\|U-V\|^\a_{\bS_p}.
$$
\end{thm}

The proofs of the above results are almost the same as in the self-adjoint case.
The only difference is that we have to use \rf{UVI} instead of the corresponding inequality for self-adjoint operators. 

We proceed now to higher order differences. Let $U$ be a unitary operator and 
$A$ a self-adjoint operator. We are going to study properties of the following higher order differences
\bay
\label{mya}
\sum_{k=0}^m(-1)^k\left(\begin{matrix}m\\k\end{matrix}\right)
f\big(e^{{\rm i}kA}U\big).
\ey
As we have already mentioned in the introduction to this section, such finite differences can be expressed as a linear combination of multiple operator integrals of multiplicity at most $m+1$. We refer the reader to \cite{AP2}, Th. 5.2. For simplicity, we state the formula in the case $m=3$. Let 
$f\in B_{\be1}^2$. Let $U_1$, $U_2$, and $U_3$ be unitary operators. Then
\begin{align}
\label{m=2}
f(U_1)&-2f(U_2)+f(U_3)\\[.2cm]
&=
2\iiint(\dg^2f)(\z,\t,\up)\,dE_1(\z)(U_1-U_2)\,dE_2(\t)(U_2-U_3)\,dE_3(\up)\nonumber\\[.2cm]
&+\iint(\dg f)(\z,\t)\,dE_1(\z)(U_1-2U_2+U_3)\,dE_3(\t).\nonumber
\end{align}
Let $U_1=U$, $U_2=e^{{\rm i}A}U$, and $U_3=e^{2{\rm i}A}U$.

\begin{lem}
\label{r3p}
Let $\fI$ be a normed ideal such that $\fI^{\{1/2\}}$ is also a normed ideal. 
If $f$ is a trigonometric polynomial of degree $d$ and $A\in\fI$, then
$f(U)-2f\big(e^{{\rm i}A}U\big)+f\big(e^{2{\rm i}A}U\big)\in\fI^{\{1/2\}}$ and
$$
\Big\|f(U)-2f\big(e^{{\rm i}A}U\big)+
f\big(e^{2{\rm i}A}U\big)\Big\|_{\fI^{\{1/2\}}}
\le\const\cdot\, d^2\,\|f\|_{L^\be}\|A\|_\fI^2.
$$
Moreover, the constant does not depend on $\fI$.
\end{lem}

\Pf Let $U_1=U$, $U_2=e^{{\rm i}A}U$, and $U_3=e^{2{\rm i}A}U$. By \rf{Bok}, we have
\begin{align*}
&\left\|
\iiint(\dg^2f)(\z,\t,\up)\,dE_1(\z)(U_1-U_2)\,dE_2(\t)(U_2-U_3)\,dE_3(\up)
\right\|_{\fI^{\{1/2\}}}\\[.2cm]
&\le\const\cdot\, d^2\,\|f\|_{L^\be}\|U_1-U_2\|_\fI\|U_2-U_3\|_\fI.
\end{align*}
Clearly,
$$
\|U_1-U_2\|_\fI=\|U_2-U_3\|_\fI=\big\|I-e^{{\rm i}A}\big\|_\fI
\le\const\|A\|_\fI.
$$
On the other hand, by \rf{Bp},
$$
\left\|\iint(\dg f)(\z,\t)\,dE_1(\z)(U_1-2U_2+U_3)\,dE_3(\t)
\right\|_{\fI^{\{1/2\}}}
\le\const\cdot\, d\,\|U_1-2U_2+U_3\|_{\fI^{\{1/2\}}}
$$
and
$$
\|U_1-2U_2+U_3\|_{\fI^{\{1/2\}}}
=\big\|(I-e^{{\rm i}A})^2\big\|_{\fI^{\{1/2\}}}\le\const\|A\|_\fI^2.
$$
The result follows now from \rf{m=2}. $\bl$

In the general case the following inequality holds:
\bay
\label{liu}
\left\|\sum_{k=0}^m(-1)^k\left(\begin{matrix}m\\k\end{matrix}\right)
f\big(e^{{\rm i}kA}U\big)\right\|_{\fI^{\{1/m\}}}
\le\const\cdot\, d^m\,\|f\|_{L^\be}\|A\|^m_\fI,
\ey
whenever $\fI$ is a normed ideal such that $\fI^{\{1/m\}}$ is also a normed ideal. This follows from an analog of formula \rf{m=2} for higher order differences, see \cite{AP2}, Th. 5.2.

We state the remaining results in this section without proofs. The proofs are practically the same as in the self-adjoint case. The only difference is that
instead of inequality \rf{li}, one has to use inequality \rf{liu}.

\begin{thm}
\label{vypu}
Let $\a>0$ and $m-1\le\a<m$. There exists a positive number $c$ such that for every $l\ge0$, $p\in[m,\be)$, $f\in\L_\a$, and for arbitrary unitary operator $U$  self-adjoint operator $A$, the following inequality holds:
$$
s_j\left(
\sum_{k=0}^m(-1)^k\left(\begin{matrix}m\\k\end{matrix}\right)
f\big(e^{{\rm i}kA}U\big)
\right)\le c\,\|f\|_{\L_\a}(1+j)^{-\a/p}\|A\|^\a_{\bS_p^l}
$$
for $j\le l$.
\end{thm}

\begin{thm}
\label{vslu}
Let $\a>0$ and $m-1\le\a<m$. There exists a positive number $c$ such that for every $f\in\L_\a$, and for an arbitrary unitary operator $U$ and an arbitrary self-adjoint operator $A$ of class $\bS_m$, the operator 
{\em\rf{mya}} belongs to
$\bS_{\frac{m}{\a},\be}$ and the following inequality holds:
$$
\left\|\sum_{k=0}^m(-1)^k\left(\begin{matrix}m\\k\end{matrix}\right)
f\big(e^{{\rm i}kA}U\big)\right\|_{\bS_{\frac{m}{\a},\be}}
\le c\,\|f\|_{\L_\a}\|A\|_{\bS_m}^\a.
$$
\end{thm}

\begin{thm}
\label{p=mu}
Let $\a>0$ and $m-1\le\a\le m$. There exists a positive number $c$ such that for every $f\in B_{\be1}^\a$, and for an arbitrary unitary operator $U$  and
an arbitrary self-adjoint operator $A$ of class $\bS_m$, the operator {\em\rf{mya}} belongs to
$\bS_{\frac{m}{\a}}$ and the following inequality holds:
$$
\left\|\sum_{k=0}^m(-1)^k\left(\begin{matrix}m\\k\end{matrix}\right)
f\big(e^{{\rm i}kA}U\big)\right\|_{\bS_{\frac{m}{\a}}}
\le c\,\|f\|_{B_{\be1}^\a}\|A\|_{\bS_m}^\a.
$$
\end{thm}

\begin{thm}
\label{vsinu}
Let $\a>0$ and $m-1\le\a<m$. There exists a positive number $c$ such that for every $f\in\L_\a$, and for arbitrary unitary operator $U$ and bounded self-adjoint operator $A$, the following inequality holds:
$$
s_j\left(\left|\sum_{k=0}^m(-1)^k\left(\begin{matrix}m\\k\end{matrix}\right)
f\big(e^{{\rm i}kA}U\big)\right|^{m/\a}\right)\le c\,
\|f\|_{\L_\a}^{m/\a}\s_j\big(|A|^m\big),\quad j\ge0.
$$
\end{thm}

\begin{thm}
\label{idu}
Let $\a>0$ and $m-1\le\a<m$. There exists a positive number $c$ such that for every $f\in\L_\a$, every quasinormed ideal $\fI$
with $\b_\fI<m^{-1}$, and for arbitrary unitary operator $U$ and bounded self-adjoint operator $A$, the following inequality holds:
$$
\left\|\,\left|\sum_{k=0}^m(-1)^k\left(\begin{matrix}m\\k\end{matrix}\right)
f\big(e^{{\rm i}kA}U\big)\right|^{1/\a}\right\|_\fI
\le c\,\bs{C}^{1/m}_{\fI^{\{1/m\}}}\|f\|_{\L_\a}^{1/\a}\|A\|_\fI.
$$
\end{thm}

\begin{thm}
\label{sjmu}
Let $\a>0$, $m-1\le\a<m$, and $m<p<\be$. There exists a positive number $c$ such that for every $f\in\L_\a$, every $l\in\Z_+$, and for arbitrary unitary operator $U$ and bounded self-adjoint operator $A$, the following inequality holds:
$$
\sum_{j=0}^l\left(s_j\left(\left|
\sum_{k=0}^m(-1)^k\left(\begin{matrix}m\\k\end{matrix}\right)
f\big(e^{{\rm i}kA}U\big)
\right|^{1/\a}\right)\right)^p
\le c\,\|f\|_{\L_\a}^{p/\a}\sum_{j=0}^l\big(s_j(A)\big)^p.
$$
\end{thm}

\begin{thm}
\label{poslu}
Let $\a>0$, $m-1\le\a<m$, and $m<p<\be$. There exists a positive number $c$ such that for every $f\in\L_\a$, for an arbitrary unitary operator $U$, and an arbitrary self-adjoint operator $A$ of class $\bS_p$, the following inequality holds:
$$
\left\|
\sum_{k=0}^m(-1)^k\left(\begin{matrix}m\\k\end{matrix}\right)
f\big(e^{{\rm i}kA}U\big)
\right\|_{\bS_{p/\a}}
\le c\,\|f\|_{\L_\a}\|A\|_{\bS_p}^\a.
$$
\end{thm}

\

\section{\bf The case of contractions}
\setcounter{equation}{0}
\label{con}

\

In this section we obtain analogs of the results of Sections \ref{sa}
and \ref{u} for contractions. To obtain desired estimates, we use multiple operator integrals with respect to semi-spectral measures. 

Suppose that $T$ and $R$ are contractions on Hilbert space and $f$ is a function in the disk-algebra $C_A$ (i.e., $f$ is analytic in $\dd$ and continuous in $\clos\dd$). We are going to study properties of differences
\bay
\label{vpr}
\sum_{k=0}^m(-1)^k\left(\begin{matrix}m\\k\end{matrix}\right)f\left(T+\frac{k}{m}(T-R)\right),\quad m\ge1.
\ey
In particular, when $m=1$, we obtain first order differences $f(T)-f(R)$. In this section we are not going to state separately results for first order differences. They can be obtained from the general results by putting $m=1$.

It was shown in \cite{AP2} that 
\begin{align}
\label{mic}
&\sum_{k=0}^m(-1)^k\left(\begin{matrix}m\\k\end{matrix}\right)f\left(T+\frac{k}{m}(T-R)\right)\\[.2cm]
&=
\frac{m!}{m^m}\underbrace{\int\cdots\int}_{m+1}
(\dg^{m}f)(\z_1,\cdots,\z_{m+1})
\,d\E_1(\z_1)(T-R)\cdots(T-R)\,d\E_{m+1}(\z_{m+1}),\nonumber
\end{align}
where $\E_k$ is a semi-spectral measure of $T+\frac km(T-R)$.

Suppose now that $\fI$ is a normed ideal such that $\fI^{\{1/m\}}$ is also a normed ideal. It follows from \rf{mic} and \rf{Boke} that for an arbitrary trigonometric polynomial $f$ of degree $d$,
\bay
\label{hdc}
\left\|
\sum_{k=0}^m(-1)^k\left(\begin{matrix}m\\k\end{matrix}\right)f\left(T+\frac{k}{m}(T-R)\right)
\right\|_{\fI^{\{1/m\}}}
\le\const\cdot\,d^m\,\|f\|_{L^\be}\|T-R\|_\fI^m,
\ey   
where the constant can depend only on $m$.

We state the results without proofs. The proofs are almost the same as in the self-adjoint case. The only difference is that to estimate higher order differences, we should use inequality \rf{hdc}.

\begin{thm}
\label{vypc}
Let $\a>0$ and $m-1\le\a<m$. There exists a positive number $c$ such that for every $l\ge0$, $p\in[m,\be)$, $f\in\big(\L_\a\big)_+$, and for arbitrary contractions $T$ and $R$ on Hilbert space, the following inequality holds:
$$
s_j\left(\sum_{k=0}^m(-1)^k\left(\begin{matrix}m\\k\end{matrix}\right)f\left(T+\frac{k}{m}(T-R)\right)\right)
\le c\,\|f\|_{\L_\a}(1+j)^{-\a/p}\|T-R\|^\a_{\bS_p^l}
$$
for $j\le l$.
\end{thm}

\begin{thm}
\label{vslc}
Let $\a>0$ and $m-1\le\a<m$. There exists a positive number $c$ such that for every $f\in\big(\L_\a\big)_+$, and for  arbitrary contractions $T$ and $R$ with 
$T-R\in\bS_m$, the operator {\em\rf{vpr}} belongs to
$\bS_{\frac{m}{\a},\be}$ and the following inequality holds:
$$
\left\|
\sum_{k=0}^m(-1)^k\left(\begin{matrix}m\\k\end{matrix}\right)f\left(T+\frac{k}{m}(T-R)\right)
\right\|_{\bS_{\frac{m}{\a},\be}}
\le c\,\|f\|_{\L_\a}\|T-R\|_{\bS_m}^\a.
$$
\end{thm}

\begin{thm}
\label{p=mc}
Let $\a>0$ and $m-1\le\a\le m$. There exists a positive number $c$ such that for every $f\in\big(B_{\be1}^\a\big)_+$, and for arbitrary contractions $T$  and
$R$ with $T-R\in\bS_m$, the operator {\em\rf{vpr}} belongs to
$\bS_{\frac{m}{\a}}$ and the following inequality holds:
$$
\left\|
\sum_{k=0}^m(-1)^k\left(\begin{matrix}m\\k\end{matrix}\right)f\left(T+\frac{k}{m}(T-R)\right)
\right\|_{\bS_{\frac{m}{\a}}}
\le c\,\|f\|_{B_{\be1}^\a}\|T-R\|_{\bS_m}^\a.
$$
\end{thm}

\begin{thm}
\label{vsinc}
Let $\a>0$ and $m-1\le\a<m$. There exists a positive number $c$ such that for every $f\in\big(\L_\a\big)_+$, and for arbitrary contractions $T$ and $R$, the following inequality holds:
$$
s_j\left(\left|
\sum_{k=0}^m(-1)^k\left(\begin{matrix}m\\k\end{matrix}\right)f\left(T+\frac{k}{m}(T-R)\right)
\right|^{m/\a}\right)\le c\,
\|f\|_{\L_\a}^{m/\a}\s_j\big(|T-R|^m\big),\quad j\ge0.
$$
\end{thm}

\begin{thm}
\label{idc}
Let $\a>0$ and $m-1\le\a<m$. There exists a positive number $c$ such that for every $f\in\big(\L_\a\big)_+$, every quasinormed ideal $\fI$
with $\b_\fI<m^{-1}$, and for arbitrary contractions $T$ and $R$, the following inequality holds:
$$
\left\|\,\,\left|
\sum_{k=0}^m(-1)^k\left(\begin{matrix}m\\k\end{matrix}\right)f\left(T+\frac{k}{m}(T-R)\right)
\right|^{1/\a}\right\|_\fI
\le c\,\bs{C}^{1/m}_{\fI^{\{1/m\}}}\|f\|_{\L_\a}^{1/\a}\|T-R\|_\fI.
$$
\end{thm}

\begin{thm}
\label{sjmc}
Let $\a>0$, $m-1\le\a<m$, and $m<p<\be$. There exists a positive number $c$ such that for every $f\in\big(\L_\a\big)_+$, every $l\in\Z_+$, and for arbitrary contractions $T$ and $R$, the following inequality holds:
$$
\sum_{j=0}^l\left(s_j\left(\left|
\sum_{k=0}^m(-1)^k\left(\begin{matrix}m\\k\end{matrix}\right)f\left(T+\frac{k}{m}(T-R)\right)
\right|^{1/\a}\right)\right)^p
\le c\,\|f\|_{\L_\a}^{p/\a}\sum_{j=0}^l\big(s_j(T-R)\big)^p.
$$
\end{thm}

\begin{thm}
\label{poslc}
Let $\a>0$, $m-1\le\a<m$, and $m<p<\be$. There exists a positive number $c$ such that for every $f\in\big(\L_\a\big)_+$, for arbitrary contractions $T$ and $R$ with $T-R\in\bS_p$, the following inequality holds:
$$
\left\|
\sum_{k=0}^m(-1)^k\left(\begin{matrix}m\\k\end{matrix}\right)f\left(T+\frac{k}{m}(T-R)\right)
\right\|_{\bS_{p/\a}}
\le c\,\|f\|_{\L_\a}\|T-R\|_{\bS_p}^\a.
$$
\end{thm}

%






\

\section{\bf Finite rank perturbations and necessary conditions. Unitary operators}
\setcounter{equation}{0}
\label{fru}

\

In this sections we study the case of finite rank perturbations of unitary operators. We also obtain some necessary conditions. In particular we show that
the assumptions that $\rank(U-V)=1$ and $f\in\L_\a$, $0<\a<1$, do not imply that $f(U)-f(V)\in\bS_{1/\a}$.

Let us introduce the notion of Hankel operators. For $\f\in L^\be(\T)$, the {\it Hankel operator} $H_\f$ from the Hardy class $H^2$ to $H^2_-\df L^2\ominus H^2$ is defined by
$$
H_\f g=\pp_-\f g,\quad g\in H^2,
$$
where $\pp_-$ is the orthogonal projection from $L^2$ onto $H^2_-$. Note that the operator $H_\f$ has Hankel matrix 
$$
\G_\f\df\{\hat \f(-j-k)\}_{j\ge1,k\ge0}
$$ 
with respect to the orthonormal bases $\{z^k\}_{k\ge0}$ and $\{\bar z^j\}_{j\ge1}$ of $H^2$ and $H^2_-$.

We need the following description of Hankel operators of class $\bS_p$ that was obtained in
\cite{Pe1} for $p\ge1$ and \cite {Pe0} and \cite{Se} for $p<1$ (see also \cite{Pe5}, Ch. 6):
\bay
\label{Han}
H_\f\in\bS_p\quad\Longleftrightarrow\quad\pp_-\f\in B_p^{1/p},\quad0<p<\be.
\ey

The following result gives us a necessary condition on $f$ for the assumption  $U-V\in\bS_1$ to imply that $f(U)-f(V)\in\bS_{1/\a}$.

\begin{thm}
\label{r1H}
 Suppose that $0<p<\infty$. Let $f$ be a continuous function on $\T$
such that $f(U)-f(V)\in \bS_p$, whenever
$U$ and $V$ are unitary operators with $\rank(U-V)=1$.
Then $f\in B_p^{1/p}$.
\end{thm}

\Pf
Consider the operators $U$ and $V$ on the space $L^2(\T)$ with respect to normalized Lebesgue measure on $\T$ defined by
$$
Uf=\bar z f\quad\mbox{and}\quad Vf=\bar zf-2(f,\1)\bar z,\quad f\in L^2.
$$
It is easy to see that both $U$ and $V$ are unitary operators and 
$$
\rank(V-U)=1.
$$
It is also easy to verify that for $n\ge0$,
$$
V^nz^j=\left\{
\begin{array}{ll}z^{j-n},&j\ge n,\\[.2cm]
-z^{j-n},&0\le j<n,\\[.2cm]
z^{j-n},&j<0.
\end{array}\right.
$$
It follows that for $f\in C(\T)$, we have
\begin{align*}
\big((f(V)-f(U))z^j,z^k\big)&=\sum_{n>0}\hat f(n)\big((V^nz^j,z^k)-(z^{j-n},z^k)\big)\\[.2cm]
&+\sum_{n<0}\hat f(n)\big((V^nz^j,z^k)-(z^{j-n},z^k)\big)\\[.2cm]
&=-2\left\{\begin{array}{ll}\hat f(j-k),&j\ge0,~k<0,\\[.2cm]
\hat f(j-k), &j<0,~k\ge0,\\[.2cm]
0,&\mbox{otherwise}.
\end{array}\right.
\end{align*}

If $f(U)-f(V)\in\bS_p$, it follows that the operators on $\ell^2$
with Hankel matrices 
$$
\{\hat f(j+k)\}_{j\ge0,k\ge1}\qm\{\hat f(-j-k)\}_{j\ge0,k\ge1}
$$
belong to $\bS_p$. It follows now from \rf{Han} that $f\in B_{p}^{1/p}$.
$\bl$

\medskip

{\bf Remark.} Recall that Theorem \ref{S1u} says that under the assumptions
$U-V\in\bS_1$ and $f\in\L_\a$, $0<\a<1$, the operator $f(U)-f(V)$ belongs to 
$\bS_{\frac1\a,\be}$. On the other hand, Theorem \ref{S1Bu} shows 
that the slightly stronger condition $f\in B_{\be1}^\a$ implies
that $f(U)-f(V)\in\bS_1$. 
However, the above theorem tells us that even under the much stronger assumption 
$\rank(U-V)=1$ the condition $f\in\L_\a$ does not imply that 
\lb$f(U)-f(V)\in\bS_1$. Indeed, $\L_\a\not\subset B_{1/\a}^\a$. This follows from 
the fact that  
\bay
\label{L1a}
\sum\limits_{k\ge0}a_kz^{2^k}\in\L_\a\quad\Longleftrightarrow\quad
\big\{2^{\a k}a_k\big\}_{k\ge0}\in\ell^\be
\ey
and from the fact that 
\bay
\label{B1a}
\sum\limits_{k\ge0}a_kz^{2^k}\in  B_{1/\a}^\a\quad\Longleftrightarrow\quad
\big\{2^{\a k}a_k\big\}_{k\ge0}\in\ell^{1/\a}.
\ey
Both \rf{L1a} and \rf{B1a} follows easily from \rf{bes}.

\medskip

Note that the proof of Theorem \ref{r1H} shows that if
$U$ and $V$ are the unitary operators constructed
in the proof of Theorem \ref{r1H} and $\frak I$ is a quasinormed ideal, then \lb$f(U)-f(V)\in\fI$
if and only if both $H_f$ and $H_{\ov f}$ belong to $\fI$.

The following result is closely related to Theorem \ref{S1u}, it shows that if we replace the assumption $U-V\in\bS_1$ with the stronger assumption
$\rank(U-V)<+\infty$, we can obtain the same conclusion for all $\a>0$.

\begin{thm}
\label{r1u}
Let $0<\a<\be$ and 
let $U$ and $V$ be unitary operators such that $\rank(U-V)<+\infty$.
Then $f(U)-f(V)\in {\boldsymbol S}_{\frac1\a,\infty}$ for
every function $f\in\L_\a(\T)$.
\end{thm}

\Pf Let $m$ be a positive integer and let $f\in\L_\a$. By Bernstein's theorem,
we can represent $f$ in the form $f=f_1+f_2$,
where $f_1$ is a trigonometric polynomial of degree at most $m$
and $\|f_2\|_{L^\be}\le\const m^{-\a}$ (this can be deduced easily from \rf{bes}). It is easy to see that
$$
U^m-V^m=\sum\limits_{j=0}^{m-1}U^j(U-V)V^{n-1-j}.
$$
Hence, 
$$
\Range\big(f_1(U)-f_1(V)\big)\subset\sum\limits_{j=-m}^m\Range\big(U^j(U-V)\big),
$$
and so 
$$
\rank\big(f_1(U)-f_1(V)\big)\le(2m+1)\rank(U-V),
$$
while  $\|f_2(U)-f_2(V)\|\le2\|f_2\|_{L^\be}\le\const m^{-\a}$.
It follows that 
$$
s_{(2m+1)\rank(U-V)}\big(f(U)-f(V)\big)\le\const m^{-\a}.\quad\bl
$$

We can compare Theorem \ref{r1u} with the following result obtained in 
\cite{Pe3}: if $0<p\le1$, and $U$ and $V$ are unitary operators such that
$U-V\in\bS_p$, then $f(U)-f(V)\in\bS_p$ for every $f\in B_{\be p}^{1/p}$.

The following result allows us to estimate the singular values of Hankel operators with symbols in $\L_\a$.

\begin{lem}
\label{sHa}
 Let $0<\a<\be$. Then there exists a positive number $c$ such that for every  
 $f\in\L_\a(\T)$, the following inequality holds:
 $$
 s_m(H_f)\le c\,\|f\|_{\L_\a}(1+m)^{-\a}.
 $$
\end{lem}

\Pf We can represent $f$ in the form $f=f_1+f_2$, where
$f_1$ is a trigonometrical polynomial of degree at most $m$
and $\|f_2\|\le\const(1+m)^{-\a}$. Then
$\rank H_{f_1}\le m$ and $\big\|H_{f_2}\big\|\le\const(1+m)^{-\a}$
which implies the result. $\bl$

The following theorem shows that Theorems \ref{r1u} and \ref{S1u} cannot be improved.

\begin{thm}
\label{stu}
Let $\a>0$. There exist unitary operators $U$ and $V$ 
and a real function $h$ in $\L_\a$ such that
$$
\rank(U-V)=1\qm
s_m(h(U)-h(V))\ge(1+m)^{-\a},\quad m\ge0.
$$
\end{thm}

\Pf Let $U$ and $V$ be the unitary operators defined in the proof of Theorem
\ref{r1H}. 

Consider the function $g$ defined by
\bay
\label{fyag}
g(\z)\df\sum\limits_{n=1}^\infty4^{-\a n}\left(\z^{\,4^n}+\ov{\z}^{\,4^n}\right),
\quad\z\in\T.
\ey
It follows easily from \rf{bes} that
$g\in\L_\a(\T)$. By Lemma \ref{sHa}, $s_m(H_g)\le\const(1+m)^{-\a}$,
$m\ge0$. Let us obtain a lower estimate for $s_m(H_g)$.

Consider the matrix $\G_g$ of the Hankel operator $H_g$ with respect to the
standard orthonormal bases:
$$
\G_g=\{\hat g(-j-k)\}_{j\ge1,k\ge0}=\{\hat g(j+k)\}_{j\ge1,k\ge0}.
$$
Let $n\ge1$. Define the $3\cdot4^{n-1}\times3\cdot4^{n-1}$ matrix $T_n$ by
$$
T_n=\big\{\hat g\big(j+k+4^{n-1}+1\big)\big\}_{0\le j,k<3\cdot4^{n-1}}.
$$
Clearly, $4^{\a n}T_n$ is an orthogonal matrix. Hence,
$\|T_n-R\|\ge4^{-\a n}$ for every $3\cdot4^{n-1}\times3\cdot4^{n-1}$ matrix
with $\rank R<3\cdot4^{n-1}$.
The matrix $T_n$ can be considered as a submatrix of $\G_g$.
Hence $\|\G_g-R\|\ge4^{-\a n}$ for every infinite matrix $R$ with $\rank R<3\cdot4^{n-1}$.
Thus, $s_j(\G_g)\ge4^{-\a n}$ for $j<3\cdot4^{n-1}$. 

To complete the proof, it suffices to take $h=cg$ for a sufficiently large number $c$. $\bl$

In \S\,\ref{u} we have obtained sufficient conditions on a function $f$ on $\T$
for the condition $U-V\in\bS_p$ to imply that $f(U)-f(V)\in\bS_q$ for certain $p$ and $q$. We are going to obtain here necessary conditions and consider other values $p$ and $q$.

We denote by ${\bf U}(\bS_p,\bS_q)$ the set of all continuous functions $f$ on 
$\T$ such that \lb$f(U)-f(V)\in\bS_q$, whenever $U$ and $V$ are
 unitary operators such that $U-V\in\bS_p$.

We also denote by ${\bf U}_{\rm c}(\bS_p,\bS_q)$ the set of all continuous functions $f$ on $\T$ such that \lb$f(U)-f(V)\in\bS_q$, whenever $U$ and $V$ are
commuting unitary operators such that $U-V\in\bS_p$.

Obviously, both ${\bf U}(\bS_p,\bS_q)$ and ${\bf U}_{\rm c}(\bS_p,\bS_q)$ contain the set of constant functions. We say that ${\bf U}(\bS_p,\bS_q)$ 
(or ${\bf U}_{\rm c}(\bS_p,\bS_q)$) is {\it trivial} if it contains no other functions.

Recall that the space ${\rm Lip}$ of {\it Lipschitz functions} on $\T$ is defined
as the space of functions $f$ such that
$$
\|f\|_\Li\df\sup_{\z\ne\t}\frac{|f(\z)-f(\t)|}{|\z-\t|}<\be.
$$

\begin{thm}
\label{Ucu}
Let $0<p,q<+\infty$. Then
$$
{\bf U}_{\rm c}(\bS_p,\bS_q)=
\left\{\begin{array}{ll}\L_{p/q},&p<q,\\[.2cm]
\Li,&p=q.
\end{array}\right.
$$
The space ${\bf U}_{\rm c}(\bS_p,\bS_q)$ is trivial if $p>q$.
\end{thm}

\Pf It is easy to see that $f\in {\bf U}_{\rm c}(\bS_p,\bS_q)$
if and only if
for every two sequences $\{\z_n\}$ and $\{\t_n\}$ in $\T$,
\bay
\label{pqu}
\sum|\z_n-\t_n|^p<\infty \quad\Longrightarrow\quad
\sum|f(\z_n)-f(\t_n)|^q<\infty.
\ey
Clearly, the condition $|f(\z)-f(\xi)|\le\const|\z-\xi|^{p/q}$  implies \rf{pqu}.

Consider the {\it modulus of continuity} $\o_f$ associated with $f$:
$$
\o_f(\d)\df\sup\{|f(x)-f(y)|:~|x-y|<\d\},\quad\d>0.
$$
Condition \rf{pqu} obviously implies that $\o_f(\d)<\be$ for some $\d>0$, and so it is finite for all $\d>0$. We have to prove that \rf{pqu} implies
that $\o_f(\d)\le\const\cdot\, t^{p/q}$.
Assume the contrary. Then there exist two sequences 
$\{\z_n\}$ and $\{\t_n\}$ in $\T$ such that  $\z_n\not=\t_n$
for all $n$,
$$
\lim_{n\to\be}|\z_n-\t_n|^p=0\qm
\lim_{n\to\be}\frac{|f(\z_n)-f(\t_n)|^q}{|\z_n-\t_n|^p}=\be.
$$
Now the result is a consequence of the following elementary fact:

If $\{\a_k\}$ and $\{\b_k\}$ are sequences of positive numbers such that
$\lim\limits_{k\to\be}\b_k=0$ and $\lim\limits_{k\to\be}\a_k\b_k^{-1}=+\infty$, then there exists
a sequence $\{n_k\}$ of nonnegative integers such that
$\sum n_k\b_k<+\infty$ and $\sum n_k\a_k=+\infty$. $\bl$ 

\begin{cor}
\label{slu}
Let $0<p,q<+\infty$. Then
$$
{\bf U}(\bS_p,\bS_q)\subset\left\{
\begin{array}{ll}
\L_{p/q},&  p<q,\\[.2cm]
 \Li,& p=q.
 \end{array}
 \right.
 $$
The space ${\bf U}(\bS_p,\bS_q)$ is trivial if $p>q$. 
\end{cor}

Recall that in \cite{PS} it was shown that if $f$ is a Lipschitz function on $\R$ and $1<p<\be$, then $\|f(A)-f(B)\|_{\bS_p}\le\const\|A-B\|_{\bS_p}$, whenever 
$A$ and $B$ are self-adjoint operators such that $A-B\in\bS_p$. Their method can also be used to prove an analog of this result for unitary operators. We are going to use this analog of the Potapov--Sukochev theorem for unitary operators in the following result.

\begin{thm}
\label{npu}
Let $1<q\le p<+\infty$.
Then
$$
{\bf U}(\bS_p,\bS_q)=\left\{\begin{array}{ll}\L_{p/q},&p<q,
\\[.2cm]
\Li,&p=q.
\end{array}\right.
$$
\end{thm}

\Pf By Corollary \ref{slu}, it suffices to show that 
$\L_{p/q}\subset {\bf U}(\bS_p,\bS_q)$
for $p<q$ and $\Li\subset {\bf U}(\bS_p,\bS_q)$ for $p=q$. 
The fact that $\L_{p/q}\subset {\bf U}(\bS_p,\bS_q)$ for $q<p$
is a consequence of Theorem \ref{spu}.
The inclusion
$\Li\subset {\bf U}(\bS_p,\bS_q)$ for $q=p$
is the analog of the Potapov--Sukochev theorem mentioned above. $\bl$

\medskip

{\bf Remark 1.} There exists a function $f$ of class $\Li$ such that
$f\not\in{\bf U}(\bS_p,\bS_q)$ for any $p>0$ and $q\in(0,1]$. Indeed, 
if $U$ and $V$ are the unitary operators constructed in the proof of Theorem \ref{r1H}, then $\rank(U-V)=1$ and $f(U)-f(V)\in\bS_1$ if and only if $f\in B_1^1$. It suffices to take a Lipschitz function $f$ that does not belong to $B_1^1$.

\medskip

{\bf Remark 2.} Let $\a>0$. There exists a function $f$ in $\L_\a$ such that
$f\not\in{\bf U}(\bS_p,\bS_q)$ for any $p>0$ and $q\in(0,1/\a]$. Indeed, it suffices to consider the unitary operators $U$ and $V$ constructed in the proof of Theorem \ref{r1H} and take a function $f\in\L_\a$ that does not belong to
$B_{1/\a}^\a$.

\begin{thm}
\label{Lpqu}
Let $0<p,q<+\infty$.
Then $\L_{p/q}\subset {\bf U}(\bS_p,\bS_q)$
if and only if $1<p<q$.
\end{thm}

\Pf If $1<p,q<+\infty$ or $p>q$, the result follows from Corollary \ref{slu} and Theorem \ref{npu}.
On the other hand, if $p\le q$ and $p\le1$, then
$\L_{p/q}\not\subset {\bf U}(\bS_p,\bS_q)$ by Remark 2. $\bl$

\begin{thm}
\label{Lipu}
Let $0<p,q<+\infty$.
Then $\Li\subset {\bf U}(\bS_p,\bS_q)$
if and only if $1<p\le q$ or $p\le1<q$.
\end{thm}

\Pf As in the proof of Theorem \ref{Lpqu},
it suffices to consider the case $p\le1$.
It was shown in \cite{NP} that 
$\Li\subset {\bf U}(\bS_1,\bS_q)\subset {\bf U}(\bS_p,\bS_q)$
if $p\le1<q$.
It remains to apply Remark 1. $\bl$

Now we are going to obtain a quantitative refinement of Corollary \ref{slu}.
Let $f\in C(\T)$. Put
$$
\O_{f,p,q}(\d)\df
\sup\big\{\|f(U)-f(V)\|_{\bS_q}:~\|U-V\|_{\bS_p}\le\d,\quad 
U,~V\quad \text {are unitary operators}\big\}.
$$


\begin{lem}
\label{lom}
Let $U_1$ and $U_2$ be a unitary operators with
$U_1-U_2\in \bS_p$. Then there exists a unitary operator $V$
such that 
$$\|U_{1}-V\|_{\bS_p}\le\frac{\pi\|U_1-V\|_{\bS_p}}{4}
\qm 
\|U_{2}-V\|_{\bS_p}\le\frac{\pi\|U_2-V\|_{\bS_p}}{4}.
$$
\end{lem}

\Pf Clearly, there exists a self-adjoint operator $A$ such that 
$\exp({\rm i}A)=U_1^{-1}U_2$ and $\|A\|\le\pi$.
Note that $\pi|e^{{\rm i}\theta}-1|\ge2|\theta|$ for $|\theta|\le\pi$.
Hence, $\|A\|_{\bS^p}\le\frac\pi2\|U_1-U_2\|_{\bS_p}$.
It remains to put $V=U_1\exp\big(\frac{\rm i}2A\big)$. $\bl$

\begin{cor}
\label{corlom}
Let $0<q<+\infty$. Then there exists a positive number $c_q$
such that for every $p\in(0,\be)$,
$$
\O_{f,p,q}(2\d)\le c_q\,\O_{f,p,q}(\d),\quad\d>0.
$$
\end{cor}

\begin{lem}
\label{npq}
Let $0<p,q<\be$ and let $f\in{\bf U}(\bS_p,\bS_q)$.
Then
$$
\O_{f,p,q}(n^{1/p}\,\d)\ge n^{1/q}\O_{f,p,q}(\d)
$$
for every positive integer $n$.
\end{lem}

\Pf The result is trivial if $\O_{f,p,q}(\d)=0$ or $\O_{f,p,q}(\d)=\infty$.
Suppose now that $0<\O_{f,p,q}(\d)<\infty$. Fix $\e\in(0,1)$. Let $U$ and $V$ be unitary operators
such that $\|U-V\|_{\bS^p}\le\d$ and
$\|f(U)-f(V)\|_{\bS_q}\ge(1-\e)\O_{f,p,q}(\d)$.
Put $\mathcal U=\bigoplus\limits_{j=1}^nU$ and $\mathcal V=\bigoplus\limits_{j=1}^nV$
(the orthogonal sum of $n$ copies of $U$ and $V$).
Clearly, $\|\mathcal U-\mathcal V\|_{\bS_p}\le\d n^{1/p}$, and
$$
\|f(\mathcal U)-f(\mathcal V)\|_{\bS_q}\ge(1-\e)n^{1/q}\O_{f,p,q}(\d)
$$
Hence, $\O_{f,p,q}(n^{1/p}\,\d)\ge (1-\e)n^{1/q}\O_{f,p,q}(\d)$ for every $\e\in(0,1)$. $\bl$

\begin{thm}
\label{mnpq}
Let $0<p,q<\be$ and
let $f\in{\bf U}(\bS_p,\bS_q)$. Then
$\O_{f,p,q}(\d)<+\infty$ for all $\d>0$ and
$$
\lim_{\d\to0}\frac{\O_{f,p,q}(\d)}{\d^{p/q}}=\inf_{\d>0}\frac{\O_{f,p,q}(\d)}{\d^{p/q}}
\le\sup_{\d>0}\frac{\O_{f,p,q}(\d)}{\d^{p/q}}=\lim_{\d\to\infty}\frac{\O_{f,p,q}(\d)}{\d^{p/q}},
$$
where both limits exist in $[0,+\infty]$.
In particular, if $f$ is a nonconstant function, then 
$\O_{f,p,q}(\d)\le c_1\,\d^{p/q}$
for every $\d\in(0,1]$ and $\O_{f,p,q}(\d)\ge c_2\,\d^{p/q}$
for every $\d\in[1,\infty)$, where $c_1$ and $c_2$ are positive numbers.
\end{thm}

\Pf Since $\O_{f,p,q}$ is nondecreasing, Corollary \ref{corlom} implies that either $\O_{f,p,q}(\d)$ is finite
for all $\d>0$ or $\O_{f,p,q}(\d)=\infty$ for all $\d>0$.  The latter is impossible
because we would be able to find sequences of unitary operators
$\{U_j\}$ and $\{V_j\}$ such that
$$
\bigoplus_j(U_j-V_j)\in\bS_p,\quad\mbox{but}\quad
\bigoplus_j\big(f(U_j)-f(V_j)\big)\not\in\bS_q.
$$
Hence, $\O_{f,p,q}(\d)<+\infty$ for all $\d>0$. We can find a sequence 
$\{\d_j\}_{j=1}^\infty$
of positive numbers such that $\d_j\to0$ and
$\lim\limits_{j\to\infty}\d_j^{-p/q}\O_{f,p,q}(\d_j)
=\limsup\limits_{\d\to0}\d^{-p/q}\O_{f,p,q}(\d)\df a$. Fix $\e\in(0,1)$.
Then there exists a positive integer $N$ such that 
$\d_j^{-p/q}\O_{f,p,q}(\d_j)\ge(1-\e)a$
for all $j>N$. Lemma \ref{npq} implies 
$\O_{f,p,q}(n^{1/p}\,\d_j)\ge(1-\e)a(n^{1/p}\,\d_j)^{p/q}$
for all $j>N$ and $n>0$. Hence, $\O_{f,p,q}(\d)\ge(1-\e)a\d^{p/q}$ for all
$\d>0$ and $\e\in(0,1)$. Thus $\O_{f,p,q}(\d)\ge a\d^{p/q}$ for all
$\d>0$ and $\lim\limits_{\d\to0}\d^{-p/q}\O_{f,p,q}(\d)=a$.
In the same way we can prove that
$$
\sup_{\d>0}\frac{\O_{f,p,q}(\d)}{\d^{p/q}}=\lim_{\d\to0}\frac{\O_{f,p,q}(\d)}{\d^{p/q}}.\quad\bl
$$

\

\section{\bf Finite rank perturbations and necessary conditions.\\ Self-adjoint operators}
\setcounter{equation}{0}
\label{frsa}

\

\newcommand\vL{\bs{\varLambda}}

We are going to obtain in this section analogs of the results of the previous section in the case of self-adjoint operators. We obtain estimates for 
$f(A)-f(B)$ in the case when $\rank(A-B)<\be$. We also obtain some necessary  conditions. In particular, we show that $f(A)-f(B)$ does not have to belong to $\bS_{1/\a}$ under the assumptions 
$\rank(A-B)=1$ and $f\in\L_\a(\R)$. 

However, there is a distinction between the case of unitary operators and the case of self-adjoint operators. To describe the class of functions $f$ on $\R$, for which $f(A)-f(B)\in\bS_q$, whenever $A-B\in\bS_p$, we have to introduce the space $\vL_\a$ of functions on $\R$ that satisfy the H\"older condition of order 
$\a$ uniformly on all intervals of length 1.

We are going to deal in this section with Hankel operators on the Hardy class
$H^2(\C_+)$ of functions analytic in the upper half-plane $\C_+$. Recall that the space $L^2(\R)$ can be represented as $L^2(\R)=H^2(\C_+)\oplus H^2(\C_-)$, where
$H^2(\C_-)$ is the Hardy class of functions analytic in the lower half-plane 
$\C_-$. We denote by $\bs{P}_+$ and $\bs{P}_-$ the orthogonal projections onto
$H^2(\C_+)$ and $H^2(\C_-)$. For a function $\f$ in $L^\be(\R)$, the Hankel operator \lb${\mathcal H}_\f:H^2(\C_+)\to H^2(\C_-)$ is defined by
$$
{\mathcal H}_\f g\df\bs{P}_-\f g,\quad g\in H^2(\C_+).
$$
As in the case of Hankel operators on the Hardy class $H^2$ of functions analytic in $\dd$, the Hankel operators ${\mathcal H}_\f$ of class $\bS_p$ can be described in terms of Besov spaces:
\bay
\label{Hanp}
{\mathcal H}_\f\in\bS_p\quad\Longleftrightarrow\quad
{\mathcal P}_-\f\in B_p^{1/p}(\R),\quad0<p<\be,
\ey
where the operator ${\mathcal P}_-$ on $L^\be(\R)$ is defined by
$$
{\mathcal P}_-\f\df\big(\pp_-(\f\circ\o)\big)\circ\o^{-1},\quad\f\in L^\be(\R),
$$
and $\o(\z)\df{\rm i}(1+\z)(1-\z)^{-1}$, $\z\in\T$. This was proved in \cite{Pe1} for $p\ge1$, and in \cite{Pe2} and \cite{Se} for $0<p<1$, see also
\cite{Pe5}, Ch. 6.

Note also that by Kronecker's theorem, ${\mathcal H}_\f$ has finite rank 
if and only if ${\mathcal P}_-\f$ is a rational function (see \cite{Pe5}, Ch. 1).

Recall that the Hilbert transform $\bs{H}$ is defined on $L^2(\R)$ by 
$\bs{H}g=-{\rm i}g_++{\rm i}g_-$, where we use the notation
$g_+\df\bs{P}_+g$ and $g_-\df\bs{P}_-g$.

\begin{thm}
\label{krsa} 
Let $A$ and $B$ be bounded self-adjoint operators on Hilbert space such that 
$\rank(A-B)<\infty$.
Then $f(A)-f(B)\in {\bS}_{\frac{1}{\a},\be}$ for
every function $f$ in $\L_\a(\R)$.
\end{thm}

\Pf Consider the Cayley transforms of $A$ and $B$:
$$
U=(A-{\rm i}I)(A+{\rm i}I)^{-1}\qm V=(B-{\rm i}I)(B+{\rm i}I)^{-1}.
$$
It is well known that $U$ and $V$ are unitary operators. Moreover, it is easy to see that $\rank(U-V)<\be$. Indeed,
\begin{align*}
(A-{\rm i}I)(A+{\rm i}I)^{-1}-(B-{\rm i}I)(B+{\rm i}I)^{-1}&=
2{\rm i}\Big((B+{\rm i}I)^{-1}-(A+{\rm i}I)^{-1}\Big)\\[.2cm]
&=2{\rm i}(A+{\rm i}I)^{-1}(A-B)(B+{\rm i}I)^{-1},
\end{align*}
and so $\rank(U-V)\le\rank(A-B)$.

Without loss of generality, we may assume that $f$ has compact support. Otherwise, we can multiply $f$ by an infinitely smooth function with compact support that is equal to 1 on an interval containing the spectra of $A$ and $B$.
Consider the function $h$ on $\T$ defined by 
$h(\z)=f\big(-{\rm i}(\z+{\rm i})(\z-{\rm i})^{-1}\big)$. Obviously, $h\in\L_\a$.

By Theorem \ref{r1H}, $h(U)-h(V)\in{\bS}_{\frac{1}{\a},\be}$. It remains to observe that $h(U)=f(A)$ and $h(V)=f(B)$. $\bl$

In \S\,\ref{sa} we have proved Theorem \ref{S1} that says that the condition 
$f\in\L_\a(\R)$
implies that $f(A)-f(B)\in\bS_{\frac1\a,\be}$, whenever $A-B\in\bS_1$. On the other hand, by Theorem \ref{S1s}, the stronger condition $f\in B_{\be1}^\a(\R)$
implies that $f(A)-f(B)\in\bS_{1/\a}$, whenever $A-B\in\bS_1$. The following result gives a necessary condition on $f$ for the assumption $A-B\in\bS_1$ to imply that $f(A)-f(B)\in\bS_{1/\a}$. It
shows that the condition $f\in\L_\a(\R)$ does not ensure
that $f(A)-f(B)\in\bS_{1/\a}$ even under the much stronger assumption that
$A-B$ has finite rank.

\begin{thm}
\label{nusa}
Let $f$ be a continuous function on $\R$ and let $p>0$. 
Suppose that $f(A)-f(B)\in\bS_p$, whenever $A$ and $B$ are
bounded self-adjoint
operators such that $\rank(A-B)<\infty$. Then 
$f\circ h \in B_p^{1/p}(\R)$
for every rational function $h$ that is real on $\R$ and has no pole at $\be$.\end{thm}

\Pf Let $\f\in L^\be(\R)$ and let $M_\f$ denote multiplication by $\var$.
For $g\in L^2(\R)$, we have
\begin{align}
\label{rat}
M_\f g-\bs{H}^{-1}M_\f \bs{H}g&=\f g+\bs{H}(\f (\bs{H}g))\\[.2cm]
&=\f -(\f g_+)_++(\f g_-)_++(\f g_+)_--(\f g_-)_-\nonumber\\[.2cm]
&=\f(g_++g_-) -(\f g_+)_++(\f g_-)_++(\f g_+)_--(\f g_-)_-\nonumber\\[.2cm]
&=2(\f g_+)_-+2(\f g_-)_+=2{\mathcal H}_\f g_++2{\mathcal H}^*_\f g_-.\nonumber
\end{align}
Hence, by \rf{Hanp}, $M_\f-\bs{H}^{-1}M_\f \bs{H}\in\bS_p$
if and only if $\f\in B_p^{1/p}(\R)$.
Moreover, by Kronecker's theorem, $\rank(M_\f-\bs{H}^{-1}\bs{H}_\f H)<+\infty$ if and only if $\f$ is a rational function.

Suppose now that $h$ is a rational function that takes real values on $\R$ and has no pole at $\be$. Define the bounded self-adjoint operators $A$ and $B$ by
$$
A\df M_h,\qm B\df\bs{H}^{-1}M_h \bs{H}.
$$
By \rf{rat}, $\rank(A-B)<\be$. Again, by \rf{rat} with $\f=f\circ h$, the
$$
f(A)-f(B)=M_{f\circ h}-\bs{H}^{-1}M_{f\circ h}\bs{H}
$$
belongs to $\bS_p$ if and only if $f\circ h \in B_p^{1/p}(\R)$. $\bl$

Note that the conclusion of Theorem \ref{nusa} implies that $f$ {\it belongs to
$B_p^{1/p}(\R)$ locally}, i.e., the restriction of $f$ to an arbitrary finite interval can be extended to a function of class $B_p^{1/p}(\R)$.

Now we are going to show that Theorem \ref{S1} cannot be improved even under the assumption that $\rank(A-B)=1$.

Denote by $L^2_{\rm e}(\R)$ the set of even functions in $L^2(\R)$
and by $L^2_{\rm o}(\R)$ the set of odd functions in $L^2(\R)$.
Clearly, $L^2(\R)=L^2_{\rm e}(\R)\oplus L^2_{\rm o}(\R)$.
Let $\f$ be an even function in $L^\be(\R)$. 
Then $L^2_{\rm e}(\R)$ and $L^2_{\rm o}(\R)$ are invariant
subspaces of the operators $M_\var$ and $H^{-1}M_\var H$.
The orthogonal projections $P_{\rm e}$ and $P_{\rm e}$ onto 
$L^2_{\rm e}(\R)$ and $L^2_{\rm o}(\R)$ are given by
$$
(P_{\rm e} g)(x)=\frac12\big(g(x)+g(-x)\big)\qm
(P_{\rm o} g)(x)=\frac12\big(g(x)-g(-x)\big).
$$

\begin{lem}
\label{pGf}
Let $\f(x)=(x^2+1)^{-1}$, $x\in\R$. Then 
$\big(\bs{H}\f\big)(x)=x(x^2+1)^{-1}$ and
$$
M_\f f-\bs{H}^{-1}M_\f\bs{H}f=\frac1\pi(f,\f)\f-\frac1\pi(f,\bs{H}\f)\bs{H}\f.
$$
In particular,
$$
M_\f f-\bs{H}^{-1}M_\f \bs{H}f
=\left\{\begin{array}{ll}\frac1\pi(f,\f)\f,&
f\,\,\,\,\text {is even},\\[.2cm]
-\frac1\pi(f,\bs{H}\f)\bs{H}\f,&
f\,\,\,\,\text {is odd}.
\end{array}\right.
$$
\end{lem}

\Pf It is easy to see that $\f_+(x)=\dfrac{\rm i}{2(x+{\rm i})}$, 
$\f_-(x)=-\dfrac{\rm i}{2(x-{\rm i})}$,
and $\big(\bs{H}\f\big)(x)=x(x^2+1)^{-1}$, $x\in\R$.
Hence,
\begin{align*}
M_\f f-\bs{H}^{-1}M_\f \bs{H}f&=
\f f+\bs{H}\f\bs{H}f=2(\f f_+)_-+2(\f f_-)_+\\[.2cm]
&=
2(\f_- f_+)_-+2(\f_+ f_-)_+\\[.2cm]
&=-{\rm i}\left(\frac{f_+}{x-{\rm i}}\right)_-+
{\rm i}\left(\frac{f_-}{x+{\rm i}}\right)_+=
-{\rm i}\frac{f_+({\rm i})}{x-{\rm i}}+
{\rm i}\frac{f_-(-{\rm i})}{x+{\rm i}}\\[.2cm]
&=-\frac1{2\pi}\Big(f,\frac1{x+{\rm i}}\Big)\frac1{x-{\rm i}}
-\frac1{2\pi}\Big(f,\frac1{x-{\rm i}}\Big)\frac1{x+{\rm i}}\\[.2cm]
&=\frac2\pi(f,\f_+)\f_-
+\frac2\pi(f,\f_-)\f_+\\[.2cm]
&=\frac1{2\pi}(f,\f+{\rm i}\bs{H}\f)(\f-{\rm i}\bs{H}\f)
+\frac1{2\pi}(f,\f-{\rm i}\bs{H}\f)(\f+{\rm i}\bs{H}\f)\\[.2cm]
&=\frac1\pi(f,\f)\f-\frac1\pi(f,\bs{H}\f)\bs{H}\f. \quad\bl
\end{align*}

\begin{cor}
\label{rHsa}
Let $\var(x)=(x^2+1)^{-1}$, $x\in\R$. Then 
$\rank(M_\f-\bs{H}^{-1}M_\f \bs{H})=2$,
$\rank(P_{\rm e}(M_\f-\bs{H}^{-1}M_\f \bs{H})P_{\rm e})=1$ and 
$\rank(P_{\rm o}(M_\f-\bs{H}^{-1}M_\f \bs{H})P_{\rm o})=1$. 
\end{cor}

\begin{lem}
\label{r1Hsa}
Let $\f$ be an even function in $L^\infty(\R)$.
Then
$$
s_n((M_\f-\bs{H}^{-1}M_\f \bs{H})P_{\rm e})\ge \sqrt2 s_n({\mathcal H}_\f)
$$
and
$$
s_n((M_\f-\bs{H}^{-1}M_\f \bs{H})P_{\rm o})\ge \sqrt2 s_n({\mathcal H}_\f).
$$
\end{lem}

\Pf Note that 
$\bs{P}_-(M_\f-\bs{H}^{-1}M_\f \bs{H})\big|H^2(\C_+)=2{\mathcal H}_\f$.
It remains to observe that $\sqrt2\bs{P}_+$ acts isometrically from
$L^2_{\rm e}(\R)$ onto $H^2(\C_+)$ and from $L^2_{\rm o}(\R)$ onto 
$H^2(\C_+)$. $\bl$

\begin{lem}
\label{ro}
There exists a function $\rho\in C^\infty(\T)$
such that $\rho(\z)+\rho({\rm i}\z)=1$, 
\lb$\rho(\z)=\rho\big(\,\ov\z\,\big)$ for all 
$\z\in\T$, and $\rho$
vanishes in a neighborhood of the set $\{-1,1\}$.
\end{lem}

\Pf Fix a function $\psi\in C^\infty(\T)$ such that $\psi$
vanishes in a neighborhood of the set $\{-1,1\}$, $\psi\ge0$, and
$\psi(\z)+\psi({\rm i}\z)>0$ for all $\z\in\T$.
Put
$$
\rho_0(\z)=
\frac{\psi(\z)+\psi(-\z)}{\psi(\z)+\psi({\rm i}\z)+\psi(-\z)+\psi(-{\rm i}\z)}.
$$
Clearly, $\rho_0$ vanishes in a neighborhood of the set $\{-1,1\}$, $\rho_0\ge0$, and
$\rho_0(\z)+\rho_0({\rm i}\z)=1$ for all $\z\in\T$.
It remains to put $\rho(\z)=\frac12\big(\rho_0(\z)+\rho_0\big(\,\ov\z\,\big)\big)$.
$\bl$

In what follows we fix such a function $\rho$.

\begin{lem} 
\label{gsn}
Let $g$ be a function in $\L_\a$ such that 
$g({\rm i}\z)=g(\z)$ for all $\z\in\T$.
Suppose that $\inf\limits_{n\ge0}(n+1)^\a s_n(H_g)>0$.
Then $\inf\limits_{n\ge0}(n+1)^\a s_n(H_{\rho g})>0$.
\end{lem}

\Pf Fix a positive $p$ such that $\a p<1$. 
Clearly, there exists a positive number $c_1$ such that
$\|H_g\|_{{\boldsymbol S}_p^n}^p\ge c_1n^{1-\a p}$ for all $n\ge0$.
Note that
$\|H_{\rho(z)g(z)}\|_{{\boldsymbol S}_p^n}=
\|H_{\rho({\rm i}z)g(z)}\|_{{\boldsymbol S}_p^n}$.
Hence, $\|H_{\rho g}\|_{{\boldsymbol S}_p^n}^p\ge \frac12c_1n^{1-\a p}$ for all $n\ge0$.
By Lemma \ref{sHa}, there exists a positive number $c_2$ such that
$\|H_{\rho g}\|_{{\boldsymbol S}_p^n}^p\le c_2n^{1-\a p}$ for all $n\ge1$.
Hence, there exists an integer $M$ such that
$\|H_{\rho g}\|_{{\boldsymbol S}_p^{Mn}}^p-\|H_{\rho g}\|_{{\boldsymbol S}_p^n}^p\ge n^{1-\a p}$
for all $n\ge1$. Note that
$$
\|H_{\rho g}\|_{{\boldsymbol S}_p^{Mn}}^p-\|H_{\rho g}\|_{{\boldsymbol S}_p^n}^p\le(M-1)n (s_n(H_{\rho g}))^p.
$$
Thus $\big(s_n(H_{\rho g})\big)^p\ge\frac{1}{M-1}n^{-\a p}$ for all $n\ge1$.
$\bl$

\begin{lem}
\label{Hga}
 There exists a real function $g_0\in\L_\a$ that
vanishes in a neighborhood of the set $\{-1,1\}$ and such that 
$g_0(\z)=g_0\big(\,\ov\z\,\big)$, $\z\in\T$, and $s_n(H_{g_0})\ge (n+1)^{-\a}$ for all $n\ge0$.
\end{lem}

\Pf Let $g$ is the function given by \rf{fyag}.
We can put $g_0\df C\rho g$ for a sufficiently large number $C$. $\bl$

\begin{thm}
\label{snas}
Let $\a>0$. Let $\f(x)=(x^2+1)^{-1}$. Consider the operators
$A$ and $B$ on $L^2_{\rm e}(\R)$ defined by 
$Ag=\bs{H}^{-1}M_\f\bs{H}g$ and $Bg=\f g$.  Then

{\em(i)} $\rank(A-B)=1$,

{\em(ii)} there exists a real bounded function $h\in\L_\a(\R)$ such that
$$
s_n\big(f(A)-f(B)\big)\ge(n+1)^{-\a},\quad n\ge0.
$$
\end{thm}

\Pf The equality $\rank(A-B)=1$ is a consequence of Corollary \ref{rHsa}.
Let $g_0$ be the function obtained in Lemma \ref{Hga}.
It is easy to see that there exists a real bounded function $h\in\L_\a(\R)$
such that $h(\var(x))=g_0\left(\dfrac{x-{\rm i}}{x+{\rm i}}\right)$.
It is well known (see \cite{Pe5}, Ch. 1, Sec. 8) that ${\mathcal H}_{h\circ\f}$
can be obtained from $H_{g_0}$ by multiplying on the left and on the right by unitary operators.
Hence, by Lemma \ref{r1Hsa},
$$
s_n(f(B)-f(A))\ge\sqrt2 s_n({\mathcal H}_{h\circ\f})=\sqrt2 s_n(H_{g_0})\ge\sqrt2(n+1)^{-\a}.\quad\bl
$$

\medskip

{\bf Remark.} The same result holds if we consider operators
$A$ and $B$ on $L^2_{\rm o}(\R)$ defined in the same way.

\medskip

In \S\,\ref{sa} we have obtained sufficient conditions on a function $f$ on $\R$
for the condition $A-B\in\bS_p$ to imply that $f(A)-f(B)\in\bS_q$ for certain $p$ and $q$. We are going to obtain here necessary conditions and consider other values $p$ and $q$.

As in the case of functions on $\T$, we consider the space $\Li(\R)$ of Lipschitz functions on $\R$ such that
$$
\|f\|_{\Li(\R)}\df\sup\left\{\frac{|f(x)-f(y)|}{|x-y|}:x,y\in\R,x\not=y\right\}<+\infty.
$$

For $\a\in(0,1]$, we denote by $\vL_\a$ the set of all functions defined on $\R$
and satisfying the H\"older condition of the order $\a$
uniformly on all intervals of a fixed length:
$$
\|f\|_{\vL_\a}\df
\sup\left\{\frac{|f(x)-f(y)|}{|x-y|^\a}:~x,y\in\R,~x\ne y,~|x-y|\le1\right\}<+\infty.
$$
Clearly, $f\in\vL_\a$ if and only if $\o_f(\d)\le\const\d^\a$ for $\d\in(0,1]$. Note that 
$$
\vL_1=\Li(\R).
$$

\begin{lem}
\label{Lipa}
Let $0<\a<1$. Then $\vL_\a=\L_\a(\R)+\Li(\R)$.
\end{lem}

\Pf The inclusion $\L_\a(\R)+\Li(\R)\subset\vL_\a$ is evident. 
Let $f\in\vL_\a$.
We can consider the piecewise linear function
$f_0$ such that
$f_0(n)=f(n)$ and $f\Big|[n,n+1]$ is linear for all $n\in\Z$.
Clearly, $f_0\in\Li(\R)$ and $f-f_0\in\L_\a(\R)$. $\bl$

Denote by ${\bf SA}(\bS_p,\bS_q)$ the set of all 
continuous functions $f$ on $\R$ such that \lb$f(B)-f(A)\in\bS_q$, whenever
$A$ and $B$ are self-adjoint  operators such that $B-A\in\bS_p$.

We also denote by ${\bf SA}_{\rm c}(\bS_p,\bS_q)$ the set of all 
continuous functions $f$ on $\R$ such that $f(B)-f(A)\in\bS_q$, whenever
$A$ and $B$ are commuting self-adjoint  operators such that $B-A\in\bS_p$.

\begin{thm}
\label{SApq}
Let $0<p,q<+\infty$. Then 
${\bf SA}_{\rm c}(\bS_p,\bS_q)=\vL_{p/q}$ for $p\le q$ and
the space ${\bf SA}_{\rm c}(\bS_p,\bS_q)$ is trivial for $p>q$.
\end{thm}

\Pf To prove the inclusion $\vL_{p/q}\subset {\bf SA}_{\rm c}(\bS_p,\bS_q)$, it suffices to observe that \lb$f\in {\bf SA}_{\rm c}(\bS_p,\bS_q)$
if and only if
for every two sequences $\{x_n\}$ and $\{y_n\}$ in $\R$
\bay
\label{qp}
\sum|x_n-y_n|^p<+\infty \quad\Longrightarrow\quad
\sum\big|f(x_n)-f(y_n)\big|^q<+\infty.
\ey
Condition \rf{qp} easily implies that $\o_f(\d)<+\infty$ for some $\d>0$, and
so for all $\d>0$.
To complete the proof, we have to prove that \rf{qp} implies
that $\o_f(\d)\le C\d^{p/q}$ for $\d\in(0,1]$.
This can be done in exactly the same way as in the case of unitary operators, see
the proof of Theorem \ref{Ucu}. $\bl$

The following result is an immediate consequence of Theorem \ref{SApq}.

\begin{thm}
\label{Spq}
Let $0<p,q<+\infty$. Then
${\bf SA}(\bS_p,\bS_q)\subset \vL_{p/q}$ for $p\le q$
and ${\bf SA}(\bS_p,\bS_q)$ is trivial for $p>q$.
\end{thm}

\begin{thm}
\label{SApqL}
Let $1<p\le q<+\infty$.
Then ${\bf SA}(\bS_p,\bS_q)=\vL_{p/q}$.
\end{thm}

\Pf In view of Theorem \ref{Spq}, we have to prove that 
$\vL_{p/q}\subset {\bf SA}(\bS_p,\bS_q)$.
In the case $p=q$ this was proved by Potapov and Sukochev \cite{PS}.
Suppose now that $p<q$. By Lemma \ref{Lipa}, it is sufficient to verify
that $\Li(\R)\subset {\bf SA}(\bS_p,\bS_q)$ and 
\lb$\L_{p/q}(\R)\subset {\bf SA}(\bS_p,\bS_q)$.
The first inclusion follows from the results of \cite{PS} as well as from
the results of \cite{NP}. Indeed,
$\Li(\R)\subset {\bf SA}(\bS_p,\bS_p)\subset {\bf SA}(\bS_p,\bS_q)$.
The inclusion $\L_{p/q}(\R)\subset {\bf SA}(\bS_p,\bS_q)$
follows from Theorem \ref{sp}. $\bl$

\begin{thm}
\label{pri}
If $0<p,q\le1$, then $\Li(\R)\not\subset{\bf SA}(\bS_p,\bS_q)$.
If $0<\a<1$, $0<p\le1$, and $0<q\le1/\a$, 
then $\L_\a(\R)\not\subset{\bf SA}(\bS_p,\bS_q)$.
\end{thm}

\Pf The result follows from Theorem \ref{nusa}. Indeed,
there exists a function in $\Li(\R)$ which does not belong to
$B^1_1(\R)$ locally, and for each $\a\in(0,1)$ there exists a function in
$\L_\a(\R)$ that does not belong to $B^\a_{1/a}(\R)$ locally.
$\bl$

\begin{thm}
\label{pqSA}
Let $0<q,p<+\infty$.
Then $\L_{p/q}(\R)\subset {\bf SA}(\bS_p,\bS_q)$
if and only if $1<p<q$.
\end{thm}

\Pf If $1<q,p<+\infty$ or $q<p$, the result follows Theorems 
\ref{SApqL} and \ref{Spq}.
If $p\le q$ and $p\le1$, then
$\L_{p/q}(\R)\not\subset {\bf SA}(\bS_p,\bS_q)$ by Theorem \ref{pri}. $\bl$.

\begin{thm} 
\label{LiSA}
Let $0<p,q<+\infty$.
Then $\Li(\R)\subset{\bf SA}(\bS_p,\bS_q)$
if and only if $1<p\le q$ or $p\le1<q$.
\end{thm}

\Pf In the same way as in the proof of Theorem \ref{pqSA}, we see
that it suffices to consider the case $p\le1$.
From the results of \cite{NP} or the results of \cite{PS} it follows that 
$\Li(\R)\subset {\bf SA}(\bS_1,\bS_q)\subset{\bf SA}(\bS_p,\bS_q)$
if $p\le1<q$.
The converse follows from Theorem \ref{pri}. $\bl$

Now we are going to obtain a quantitative refinement of Theorem \ref{Spq}.
Let $f\in C(\R)$. Put
$$
\O_{f,p,q}(\d)\df\sup\big\{\|f(A)-f(B)\|_{\bS_q}\!:~\|A-B\|_{\bS_p}\le\d,~ A,\,B\quad \text {are self-adjoint operators}\big\}.
$$

It is easy to see that given $q>0$, there exists a positive number $c_q$ such that $\O_{f,p,q}(2\d)\le c_q\,\O_{f,p,q}(\d)$.

\begin{thm}
\label{pSA}
Let $0<p,q<\be$ and let $f\in{\bf SA}(\bS_p,\bS_q)$. Then
$\O_{f,p,q}(\d)<+\infty$ for all $\d>0$ and
$$
\lim_{\d\to0}\frac{\O_{f,p,q}(\d)}{\d^{p/q}}=\inf_{\d>0}\frac{\O_{f,p,q}(\d)}{\d^{p/q}}
\le\sup_{\d>0}\frac{\O_{f,p,q}(\d)}{\d^{p/q}}=\lim_{\d\to+\infty}\frac{\O_{f,p,q}(\d)}{\d^{p/q}},
$$
(both limits exist and take values in  $[0,\infty]$).
In particular, if $f$ is a nonconstant function, then 
$\O_{f,p,q}(\d)\le c_1\,\d^{p/q}$
for every $\d\in(0,1]$ and $\O_{f,p,q}(\d)\ge c_2\,\d^{p/q}$
for every $\d\in[1,+\infty)$, where $c_1$ and $c_2$ are positive numbers.
\end{thm}

The proof of Theorem \ref{pSA} is the same as that of Theorem \ref{mnpq}.

\begin{thm} Let $f\in C(\R)$ and $p\in[1,+\infty)$. 
Then either $\Omega_{f,p,p}(\d)=+\infty$ for all $\d>0$ or $\Omega_{f,p,p}$ is a linear function.
\end{thm}

\Pf If $f\not\in{\bf SA}(\bS_p,\bS_p)$, then $\Omega_{f,p,p}(\d)=+\infty$ for all $\d>0$. 

Suppose now that $f\in{\bf SA}(\bS_p,\bS_p)$. 
By the analog of Lemma \ref{npq} for self-adjoint operators 
(it is easy to see that it holds for self-adjoint operators),
$\Omega_{f,p,p}(n\d)\ge n\Omega_{f,p,p}(\d)$ for all positive integer $n$.
On the other hand, clearly, $\Omega_{f,p,p}(n\d)\le n\Omega_{f,p,p}(\d)$ for all positive integer $n$.
Hence, $\Omega_{f,p,p}$ is a linear function. $\bl$






%

\

\section{\bf Spectral shift function for second order differences}
\setcounter{equation}{0}
\label{fss}

\

In this section we obtain trace formulae for second order differences in the case of self-adjoint operators and unitary operators.

By Theorem \ref{p=m}, if $A$ is a self-adjoint operator, $K$ is a self-adjoint operator of class $\bS_2$ and $f\in B_{\be1}^2(\R)$, then 
$f(A+K)-2f(A)+f(A-K)\in\bS_1$. We are going to obtain a formula for the trace of this operator. 

\begin{thm}
\label{fsa}
Let $A$ be a self-adjoint operator and $K$ a self-adjoint operator of class 
$\bS_2$. Then there exists a unique function $\varsigma\in L^1(\R)$ such that for every $f\in B_{\be1}^2(\R)$,
\bay
\label{var}
\trace\big(f(A+K)-2f(A)+f(A-K)\big)=\int_\R f''(x)\varsigma(x)\,d\m(x).
\ey
Moreover, $\varsigma(x)\ge0$, $x\in\R$.
\end{thm}

{\bf Definition.} The function $\varsigma$ satisfying \rf{var} is called the {\it second order spectral shift function associated with the pair} $(A,K)$.

\medskip

We are going to use the spectral shift function of Koplienko. Koplienko proved in \cite{Ko} that with each pair of a self-adjoint operator $A$ and a self-adjoint operator $K$ of class $\bS_2$, there exists a function $\eta\in L^1(\R)$ such that $\eta\ge\0$ and for every rational function $f$ with poles off $\R$, the following trace formula holds
\bay
\label{Ktf}
\trace\left(f(A+K)-f(A)-\frac{d}{dt}f(A+tK)\Big|_{t=0}\right)
=\int_\R f''(x)\eta(x)\,d\m(x).
\ey
The function $\eta$ is called the {\it Koplienko spectral shift function associated with the pair} $(A,K)$. Note that later in \cite{Pe6} it was proved that trace formula \rf{Ktf} holds for \lb$f\in B_{\be1}^2(\R)$.
Note that the derivative
$$
\frac{d}{dt}f(A+tK)\Big|_{t=0}
$$
exists under the condition $f\in B_{\be1}^1(\R)$ (see \cite{Pe4} and \cite{Pe7})
and does not have to exist under the condition $f\in B_{\be1}^2(\R)$. However,
in the case $f\in B_{\be1}^2(\R)$, by 
$$
f(A+K)-f(A)-\frac{d}{dt}f(A+tK)\Big|_{t=0}
$$
we can understand 
$$
\sum_{n\in\Z}\left(f_n(A+K)-f_n(A)-\frac{d}{dt}f_n(A+tK)\Big|_{t=0}\right),
$$
and the series converges absolutely, see \cite{Pe6}. Here, as usual, 
$f_n\df f*f_n+f*W_n^\sharp$.

\medskip

{\bf Proof of Theorem \ref{fsa}.} Let $\eta_1$ be the spectral shift function associated with the pair $(A,K)$ and let $\eta_2$ be the spectral shift function associated with the pair $(A,-K)$. We have
$$
\trace\left(f(A+K)-f(A)-\frac{d}{dt}f(A+tK)\Big|_{t=0}\right)
=\int_\R f''(x)\eta_1(x)\,d\m(x)
$$
and
$$
\trace\left(f(A-K)-f(A)-\frac{d}{dt}f(A-tK)\Big|_{t=0}\right)
=\int_\R f''(x)\eta_2(x)\,d\m(x)
$$
for $f\in B_{\be1}^2(\R)$. Taking the sum, we obtain
$$
\trace\left(f(A+K)-2f(A)+f(A-K)\right)=
\int_\R f''(x)\big(\eta_1(x)+\eta_2(x)\big)\,d\m(x).
$$
It remains to put $\varsigma\df\eta_1+\eta_2$.

Uniqueness is obvious. $\bl$

We proceed now to the case of unitary operators. Suppose that $U$ is a 
unitary operator and $\V$ is a unitary operator such that $I-\V\in\bS_2$. It follows from Theorem \ref{p=mu} that if $f\in B_{\be1}^2$, then 
$f(\V U)-2f(U)+2(\V^*U)\in\bS_1$. We are going to obtain a trace formula for this operator.

\begin{thm}
\label{fsu}
Let $U$ be a unitary operator and let $\V$ be a unitary operator such that 
$I-\V\in\bS_2$. Then there exists an integrable function $\varsigma$ on $\T$ such that 
\bay
\label{varu}
\trace\big(f(\V U)-2f(U)+2(\V^*U)\big)=\int_\T f''\varsigma\,d\m.
\ey
\end{thm}

It is easy to see that $\varsigma$ is determined by \rf{varu} modulo a constant.
It is called a {\it second order spectral shift function associated with the pair} $(U,\V)$.

We are going to use a trace formula of Neidhardt \cite{Ne}, which is an analog of the Koplienko trace formula for unitary operators. 

Suppose that $U$ and $V$ be unitary operators such that $U-V\in\bS_2$.
Then $V$ can be represented as $V=e^{{\rm i}A}U$, where $A$ is a self-adjoint operator of class $\bS_2$ whose spectrum $\s(A)$ is a 
subset of $(-\pi,\pi]$. It was shown in \cite{Ne} that one can associate with the pair $(U,V)$ a function $\eta$ in $L^1(\T)$ (a Neidhardt spectral shift function) such that if the second derivative $f''$ of a function $f$ has absolutely converging Fourier series, then
\bay
\label{nei}
\trace\left(f(V)-f(U)-
\frac{d}{ds}\Big(f\big(e^{{\rm i}sA}U\big)\Big)\Big|_{s=0}\right)=
\int_\T f''\eta\,d\m.
\ey
Later it was shown in \cite{Pe6} that formula \rf{nei} holds for an arbitrary function $f$ in $B_{\be1}^2$.

\medskip

{\bf Proof of Theorem \ref{fsu}.} We can represent $\V$ as 
$\V=e^{{\rm i}A}$, where $A$ is a self-adjoint operator of class $\bS_2$
such that the spectrum $\s(A)$ of $A$ is a subset of $(-\pi,\pi]$.

Let $V_1\df\V U$. Clearly, $V_1$ is a unitary operator and $U-V_1\in\bS_2$. We can represent $\V$ as 
$\V=e^{{\rm i}A}$, where $A$ is a self-adjoint operator of class $\bS_2$
such that $\s(A)\subset(-\pi,\pi]$. We have $V_1=e^{{\rm i}A}U$.
Let $V_2\df\V^*U$. Then $U-V_2\in\bS_2$
and $V_2=e^{-{\rm i}A}U$.

Let $\eta_1$ be the Neidhardt spectral shift function associated with $(U,V_1)$ and let $\eta_2$ be the Neidhardt spectral shift function associated with 
$(U,V_2)$. We have
$$
\trace\left(f(V_1)-f(U)-
\frac{d}{ds}\Big(f\big(e^{{\rm i}sA}U\big)\Big)\Big|_{s=0}\right)=
\int_\T f''\eta_1\,d\m
$$
and
$$
\trace\left(f(V_2)-f(U)-
\frac{d}{ds}\Big(f\big(e^{-{\rm i}sA}U\big)\Big)\Big|_{s=0}\right)=
\int_\T f''\eta_2\,d\m
$$
for $f\in B_{\be1}^2$. Taking the sum, we obtain
\begin{align*}
\trace\big(f(\V U)-2f(U)+2(\V^*U)\big)&=
\trace\big(f(V_1)-2f(U)+f(V_2)\big)
\\[.2cm]
&=\int_\T f''(\eta_1+\eta_2)\,d\m.
\end{align*}
It remains to put $\varsigma\df\eta_1+\eta_2$. $\bl$

\

\section{\bf Commutators and quasicommutators}
\setcounter{equation}{0}
\label{cqc}

\

In this section we obtain estimates for the norm of {\it quasicommutators} $f(A)Q-Qf(B)$ in terms of $\|AQ-QB\|$ for self-adjoint operators $A$ and $B$
and a bounded operator $Q$. We assume for simplicity that $A$ and $B$ are bounded. However, we obtain estimates that do not depend
on the norms of $A$ and $B$. In \cite{AP3} we will  consider the case of not necessarily bounded operators $A$ and $B$.
Note that in the special case $A=B$, this problem turns into the problem of estimating the norm of commutators $f(A)Q-Qf(A)$ in terms of 
$\|AQ-QA\|$. On the other hand, in the special case $Q=I$ the problem turns into the problem of estimating $\|f(A)-f(B)\|$ in terms
$\|A-B\|$. 

Similar results can be obtained for unitary operators and for contractions.

Birman and Solomyak (see \cite{BS5}) discovered the following formula
$$
f(A)Q-Qf(B)=\iint\frac{f(x)-f(y)}{x-y}\,dE_A(x)(AQ-QB)\,dE_B(y),
$$
whenever $f$ is a function, for which $\dg f$ is a Schur multiplier
of class $\fM(E_A,E_B)$ (see \S\,\ref{koi}).

\begin{thm}
\label{glq}
Let $0<\a<1$. There exists a positive number $c>0$ such that for every
$l\ge0$,  $p\in[1,\be)$,  $f\in\L_\a(\R)$, for arbitrary bounded self-adjoint operators $A$ and $B$ and an arbitrary bounded operator $Q$, the following inequality holds:
$$
s_j\big(f(A)Q-Qf(B)\big)\le 
c\,\|f\|_{\L_\a(\R)}(1+j)^{-\a/p}\|Q\|^{1-\a}\|AQ-BQ\|_{\bS_p^l}^\a
$$
for every $j\le l$.
\end{thm}

\Pf Clearly, we may assume that $Q\ne\0$.
As usual, $f_n=f*W_n+f*W_n^\sharp$, $n\in\Z$. Fix an integer $N$. We have
by \rf{BSi} and \rf{Be},
\begin{align*}
\left\|\sum_{n=-\be}^N\big(f_n(A)Q-Qf_n(B)\big)\right\|_{\bS_p^l}
&\le\sum_{n=-\be}^N\big\|f_n(A)Q-Qf_n(B)\big\|_{\bS_p^l}\\[.2cm]
&\le\const\sum_{n=-\be}^N2^n\|f_n\|_{L^\be}\|AQ-QB\|_{\bS_p^l}\\[.2cm]
&\le\const2^{N(1-\a)}\|f\|_{\L_\a(\R)}\|AQ-QB\|_{\bS_p^l}.
\end{align*}
On the other hand,
\begin{align*}
\left\|\sum_{n>N}\big(f_n(A)Q-Qf_n(B)\big)\right\|
&\le2\|Q\|\sum_{n>N}\|f_n\|_{L^\be}\\[.2cm]
&\le\const\|f\|_{\L_\a(\R)}\|Q\|\sum_{n>N}2^{-n\a}
\le\const2^{-N\a}\|f\|_{\L_\a(\R)}\|Q\|.
\end{align*}
Put
$$
X_N\df\sum_{n=-\be}^N\big(f_n(A)Q-Qf_n(B)\big)
\qm Y_N\df\sum_{n>N}\big(f_n(A)Q-Qf_n(B)\big).
$$
Clearly, for $j\le l$,
\begin{align*}
s_j\big(f(A)Q-Qf(B)\big)&\le s_j(X_N)+\|Y_N\|
\le(1+j)^{-\frac1p}\|AQ-QB\|_{\bS_p^l}+\|Y_N\|\\[.2cm]
&\le\const\|f\|_{\L_\a(\R)}\!\left((1+j)^{-\frac1p}2^{N(1-\a)}\|AQ-QB\|_{\bS_p^l}
+2^{-N\a}\|Q\|\right).
\end{align*}
To obtain the desired estimate, it suffices to choose the number $N$ so that
$$
2^{-N}<(1+j)^{-1/p}\|AQ-QB\|_{\bS_p^l}\|Q\|^{-1}\le2^{-N+1}.\quad\bl
$$

The proofs of the remaining results of this section are the same as those of 
the results of \S\,\ref{sa} for first order differences.

\begin{thm}
\label{S1q}
Let $0<\a<1$. There exists a positive number $c>0$ such that for every
$f\in\L_\a(\R)$, for arbitrary bounded self-adjoint operators $A$ and $B$ with $AQ-QB\in\bS_1$ and an arbitrary bounded operator $Q$, the operator $f(A)Q-Qf(B)$ belongs to 
$\bS_{\frac1\a,\be}$ and the following inequality holds:
$$
\big\|f(A)Q-Qf(B)\big\|_{\bS_{\frac1\a,\be}}\le c\,\|f\|_{\L_\a(\R)}
\|Q\|^{1-\a}\|AQ-BQ\|_{\bS_1}^\a.
$$
\end{thm}

\begin{thm}
\label{S1sq}
Let $0<\a\le1$. There exists a positive number $c>0$ such that for every
$f\in B_{\be1}^\a(\R)$, for arbitrary bounded self-adjoint operators $A$ and $B$ with $AQ-QB\in\bS_1$ and an arbitrary bounded operator $Q$, the operator $f(A)Q-Qf(B)$ belongs to 
$\bS_{1/\a}$ and the following inequality holds:
$$
\big\|f(A)Q-Qf(B)\big\|_{\bS_{1/\a}}\le c\,\|f\|_{B_{\be1}^\a(\R)}
\|Q\|^{1-\a}\|AQ-QB\|_{\bS_1}^\a.
$$
\end{thm}

\begin{thm}
\label{sigmaq}
Let $0<\a<1$. There exists a positive number $c>0$ such that for every
$f\in\L_\a(\R)$, for arbitrary bounded self-adjoint operators $A$ and $B$
and an arbitrary bounded operator $Q$ on Hilbert space, the following inequality holds:
$$
s_j\Big(\big|f(A)Q-Qf(B)\big|^{1/\a}\Big)\le c\,\|f\|_{\L_\a(\R)}^{1/\a}
\|Q\|^{\frac{1-\a}\a}\s_j(AQ-QB),\quad j\ge0.
$$
\end{thm}

\begin{thm}
\label{osnq}
Let $0<\a<1$. There exists a positive number $c>0$ such that for every
$f\in\L_\a(\R)$, for an arbitrary quasinormed ideal $\fI$ with $\b_\fI<1$, for arbitrary bounded self-adjoint operators $A$ and $B$  with $AQ-QB\in\fI$,
the operator $\big|f(A)Q-Qf(B)\big|^{1/\a}$ belongs to $\fI$
and the following inequality holds:
$$
\Big\|\,\big|f(A)Q-Qf(B)\big|^{1/\a}\Big\|_\fI\le c\,\bs{C}_\fI
\|f\|_{\L_\a(\R)}^{1/\a}\|Q\|^{\frac{1-\a}\a}\|AQ-QB\|_\fI.
$$
\end{thm}

\begin{thm}
\label{splq}
Let $0<\a<1$ and $1<p<\be$. There exists a positive number $c$ such that
for every $f\in\L_\a(\R)$, every $l\in\Z_+$, for arbitrary bounded self-adjoint operators $A$ and $B$ and an arbitrary bounded operator $Q$, the following inequality holds:
$$
\sum_{j=0}^l\left(s_j\Big(\big|f(A)Q-Qf(B)\big|^{1/\a}\Big)\right)^p\le
c\,\|f\|_{\L_\a(\R)}^{p/\a}\|Q\|^{p\frac{1-\a}\a}\sum_{j=0}^l\big(s_j(AQ-QB)\big)^p.
$$
\end{thm}

\begin{thm}
\label{spq}
Let $0<\a<1$ and $1<p<\be$. There exists a positive number $c$ such that
for every $f\in\L_\a(\R)$, for arbitrary bounded self-adjoint operators $A$ and $B$, and for an arbitrary bounded operator $Q$, the operator $f(A)Q-Qf(B)$ belongs to $\bS_{p/\a}$ and the following inequality holds:
$$
\big\|f(A)Q-Qf(B)\big\|_{\bS_{p/\a}}
\le c\,\|f\|_{\L_\a(\R)}\|Q\|^{1-\a}\|AQ-QB\|^\a_{\bS_p}.
$$
\end{thm}

\

\

\noindent
\begin{tabular}{p{9cm}p{15cm}}
A.B. Aleksandrov & V.V. Peller \\
St-Petersburg Branch & Department of Mathematics \\
Steklov Institute of Mathematics  & Michigan State University \\
Fontanka 27, 191023 St-Petersburg & East Lansing, Michigan 48824\\
Russia&USA
\end{tabular}
\end{document}

%% file: classstartscr.tex
\newcommand{\lb}{\linebreak}

\renewcommand{\a}{\alpha}
\renewcommand{\b}{\beta}
\newcommand{\g}{\gamma}
\renewcommand{\d}{\delta}
\newcommand{\e}{\varepsilon}

\newcommand{\z}{\zeta}

\renewcommand{\l}{\lambda}

\newcommand{\s}{\sigma}
\renewcommand{\t}{\tau}
\newcommand{\var}{\varphi}
\newcommand{\f}{\varphi}
\renewcommand{\o}{\omega}

\newcommand{\G}{\Gamma}
\newcommand{\D}{\Delta}
\renewcommand{\L}{\Lambda}

\renewcommand{\O}{\Omega}

\newcommand{\B}{{\mathscr B}}

\newcommand{\E}{{\mathscr E}}

\newcommand{\F}{{\mathscr F}}

\newcommand{\h}{{\mathscr H}}

\newcommand{\K}{{\mathscr K}}

\newcommand{\X}{{\mathscr X}}
\newcommand{\Y}{{\mathscr Y}}

\newcommand{\V}{{\mathscr V}}

\newcommand{\1}{{\bf 1}}

\newcommand{\C}{{\Bbb C}}
\newcommand{\T}{{\Bbb T}}
\newcommand{\pp}{{\Bbb P}}
\newcommand{\dd}{{\Bbb D}}
\newcommand{\R}{{\Bbb R}}
\newcommand{\Z}{{\Bbb Z}}

\newcommand{\0}{{\boldsymbol{0}}}

\newcommand{\bs}{\boldsymbol}

\newcommand{\m}{{\boldsymbol m}}
\newcommand{\bS}{{\boldsymbol S}}

\newcommand{\rf}[1]{(\ref{#1})}

\newcommand{\df}{\stackrel{\mathrm{def}}{=}}

\newcommand{\spn}{\operatorname{span}}
\newcommand{\supp}{\operatorname{supp}}
\newcommand{\clos}{\operatorname{clos}}

\newcommand{\trace}{\operatorname{trace}}
\newcommand{\rank}{\operatorname{rank}}
\newcommand{\const}{\operatorname{const}}

\newcommand{\eeq}{\end{equation}}
\newcommand{\beq}{\begin{equation}}
\newcommand{\bay}{\begin{eqnarray}}
\newcommand{\ba}{\begin{align*}}
\newcommand{\ea}{\end{align*}}
\newcommand{\ey}{\end{eqnarray}}
\newcommand{\bey}{\begin{eqnarray*}}
\newcommand{\eey}{\end{eqnarray*}}

\newcommand{\imp}{\Rightarrow}
\newcommand{\be}{\infty}

\newcommand{\bl}{\blacksquare}

\newcommand{\Range}{\operatorname{Range}}

\newcommand{\Pf}{{\bf Proof. }}

\newcommand{\ov}{\overline}

\newtheorem{thm}{\hspace{\parindent}Theorem}[section]

\newtheorem{cor}[thm]{\hspace{\parindent}Corollary}
\newtheorem{lem}[thm]{\hspace{\parindent}Lemma}